\documentclass[11pt]{elsarticle}
\usepackage{lineno, array, booktabs, multirow, tabulary, ulem}
\modulolinenumbers[5]
\newcommand{\head}[1]{\textnormal{\textbf{#1}}}
\input{cf_macros.sty}

\begin{document} 

\begin{frontmatter}

\title{Seafloor identification in sonar imagery via simulations of Helmholtz equations and discrete optimization}

\author[UT]{Bj\"orn Engquist}
\ead{engquist@math.utexas.edu}

\author[GT]{Christina Frederick\corref{mycorrespondingauthor}}
\cortext[mycorrespondingauthor]{Corresponding author}
\ead{cfrederick6@math.gatech.edu}

\author[navy]{Quyen Huynh}
\ead{quyen\_huynh@brown.edu}

\author[GT]{Haomin Zhou}
\ead{hmzhou@gatech.edu}

\address[UT]{Department of Mathematics, The University of Texas at Austin}
\address[GT]{School of Mathematics, Georgia Institute of Technology}
\address[navy]{Institute for Brain and Neural Systems, Brown University}


\begin{abstract}We present a multiscale approach for identifying features in ocean beds by solving inverse problems in high frequency seafloor acoustics. The setting is based on Sound Navigation And Ranging (SONAR) imaging used in scientific, commercial, and military applications. The forward model incorporates multiscale simulations, by coupling Helmholtz equations and geometrical optics for a wide range of spatial scales in the seafloor geometry. This allows for detailed recovery of seafloor parameters including material type. Simulated backscattered data is generated using numerical microlocal analysis techniques. In order to lower the computational cost of the large-scale simulations in the inversion process, we take advantage of a pre-computed library of representative acoustic responses from various seafloor parameterizations. \end{abstract}


\end{frontmatter}

\newcommand{\pf}{{p_{f}}}
\newcommand{\Pf}{{\hat{p}}}
\newcommand{\ps}{p_{s}}
\newcommand{\Ps}{\hat{p}_{s}}
\newcommand{\vf}{\phi}
\newcommand{\vs}{\psi}
\newcommand{\ds}{\text{ }ds}
\newcommand{\dx}{\text{ }d\x}
\renewcommand{\r}[1]{{\color{red} #1}}
\section{Introduction}
The acquisition and interpretation of high resolution imagery of the ocean beds is needed in a wide range of oceanographic applications, for example the classification of marine habitats, sediment composition, and detection of naval mines, just to name a few \cite{Chapman1999}. Acoustic systems that penetrate 
seawater, such as sidescan sonar systems, are primarily used in underwater imaging \cite{Blondel2009}. Due to recent advancements in the resolution and speed of acquisition \cite{Groen2010}, high quality sonar data is increasingly available. This has motivated much scientific effort in developing quantitative methods for extracting details about the seafloor from the measured acoustic response.

Most methods for material classification or object identification using sonar data fall into one of two categories. In the first category, statistical models of the seafloor are developed, often from empirical training datasets \cite{Reed1989,Cervenka1993, Dobeck,  Fox1985,Jackson2007}, and used in combination with image processing based methods for classification of texture, geometric and spectral features \cite{DelRioVera2009,Reed2003}. Matched field processing is a parameter estimation technique for discriminating between the desired signal and incoherence using a-priori knowledge of the environment \cite{Baggeroer1988,Presig1994}. These techniques are successful in object detection, however they often produce significantly more false detections than the number of true targets   \cite{Groen2010}. The second category involves determining parameters in physical models of underwater acoustic wave propagation by matching the data with predictions. For example, in \cite{Burnett2015, Rzhanov2012}, a library of  representative acoustic responses is formed using computer simulations of Helmholtz equations and a label is assigned to regions within the image based on the best fit in the library. In general, the theory and simulations are well-documented for models in low frequency acoustics \cite{Frisk2009, Jackson2007}. Towards the high frequency range, most models and methods are restricted to small spatial domains.

The resolution of these images is determined by many factors, including the bandwidth of  source signal, the distance between the source array and the target, and the motion estimate of source and receivers. Understanding the high frequency response from the seafloor is essential to producing high quality images and solving problems of material classification and target identification. It is therefore important to tackle the inverse problem with accurate modeling and simulation of the entire physical wave propagation and scattering process governed by Helmholtz equations. Though there are plenty of fast Helmholtz solvers available, simulations on large domains are still cost-prohibitive.  

\begin{figure}\caption{{Wave propagation in an underwater acoustic environment consisting of a water layer and sedimentary layers. In the water layer far from the seafloor, direct ray tracing formulas can be used to approximate the traveling waves. Scattering effects close to the seafloor require accurate simulations of Helmholtz equations.}}\label{fig:SFdiagram}
\begin{center}
\includegraphics[scale=.8]{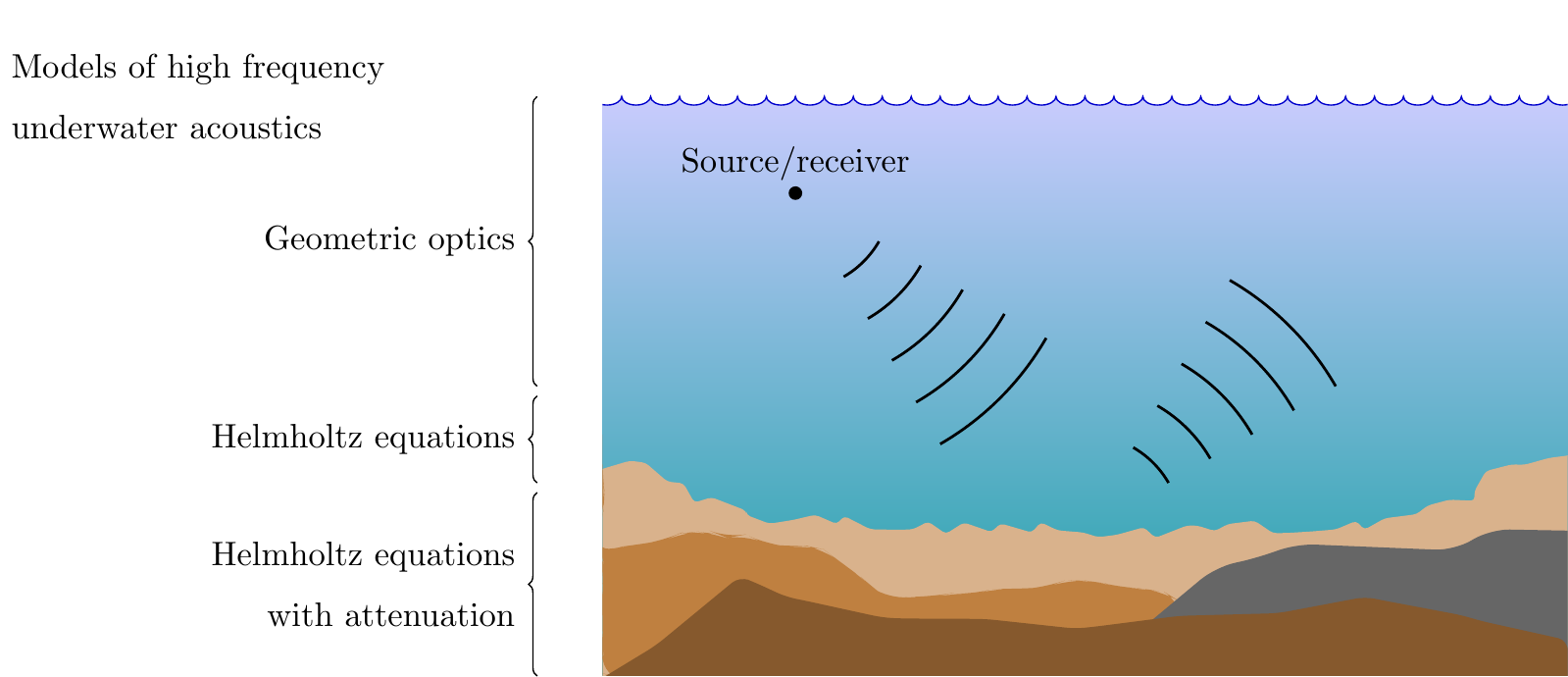}
\end{center}

\end{figure}

By taking advantage of the separation in scales between high frequency acoustic wavelengths and large spatial domains, the full scattering process can be broken down into separate stages. As depicted in Figure \ref{fig:SFdiagram}, an acoustic signal emitted from the source array travels toward the ocean bottom, generating scattered signals as a result of interactions with the seabed. The reverberations are recorded by an array of receivers near the surface and used to create sonar images. There are three main stages:  (a) the incoming wave from the source to the sea floor; (b) thin surface layer scattering; and  (c) scattering wave propagation toward the receivers. In the high frequency range, stages (a) and (c) are relatively easy to model. It should be noted that variable wave propagation velocity in the water can easily be handled by geometrical optics. For example, direct ray tracing formulas involving the attenuation rate, the distance from the source to the scattering surface,  as well as the distance between the surface and the receivers, can be used to model the pressure waves far from the seafloor \cite{Glassner1989, Jensen2011}. In \cite{Kargl2014, Kargl2012} a careful study is carried out that confirms a fast ray model for accurately simulating source-target and target-receiver propagation when the scattering field near the target is known.

{In comparison, stage (b) is much more complex and difficult to handle. Accurate simulation near the sea floor is one of the main challenges in the field because the scattering properties are influenced by a wide range of factors. For example, what are the acoustic parameters, such as wavelength and incoming angle? Does the sea floor consist primarily of clay, sand, sea grass, or rock? What is the surface geometry, especially on the scale that might produce high interference with the incoming waves?  Among all the parameters involved in the sonar model, some are given by the problem setup, for example, the acoustic wavelength and incident angle of the source ping. Other parameters are more difficult to obtain, and in fact, determining these parameters is often the main goal in many applications. Here, we aim to provide methods for identifying seabed properties, such as material type, using the recorded backscattered signal. This problem is a typical inverse problem that is often solved using many simulations of backscatter models repetitively in an iterative procedure.

One main challenge is the large computational cost associated with high frequency scattering models. Usually this is addressed by simplifying the forward models, for example, some studies assume a smooth seafloor (as shown in Figure \ref{fig:msinterface}) and approximate the scattered rays by {just following the theory of geometric optics with empirical formulas for diffuse reflections \cite{Hanrahan1993, Jensen2001, Lurton2010}}. These simplified models usually provide inadequate descriptions of the complex scattering physics caused by spatially varying material types and microtopography of the sediment layer. Also, in applications such as marine archaeology and mine detection, the targeted objects are submerged, and simple surface scattering approximations are not valid \cite{Groen2010,Kessel2002}. There are also practical challenges that influence the modeling.  For example, due to the placement of sources and receivers, typical sonar systems record data only in the backscattered direction, which is not necessarily the strongest scattering direction. The weak response is at a higher risk for noise corruption.  

This paper is provides a different methodology for material classification and object identification, taking these concerns into consideration. The relevant physics is captured by wave theory models combined with a novel technique for determining backscatter, and inversion is performed by incorporating accurate and efficient simulations.  More precisely, we assemble a library of acoustic templates representing backscatter from seafloors with varying geometric and geoacoustic parameters. Templates are created using fine scale simulations of Helmholtz equations, pre-computed offline, on partitioned small subdomains of a thin layer near the sea floor. The library is used for material classification and object identification, which are the typical inverse problems. The data is matched with predictions using optimization strategies that take advantage of multiscale wavelet representations. 

\begin{figure}\caption{On the left is a diagram depicting diffuse scattering from a rough seafloor and on the right is a diagram of reflections from a smooth seafloor. }\label{fig:msinterface}
\begin{center}
\includegraphics[height=.3\linewidth, width=.4\linewidth ]{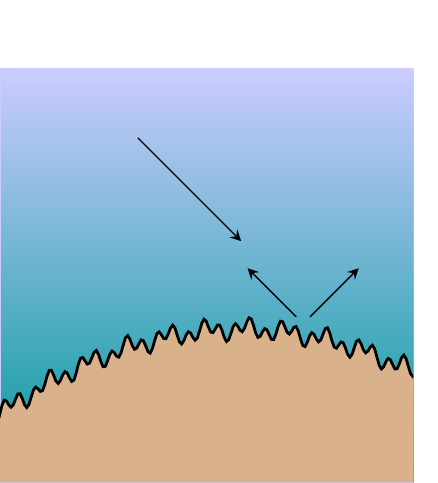} \hspace{2em}\includegraphics[height=.3\linewidth, width=.4\linewidth ]{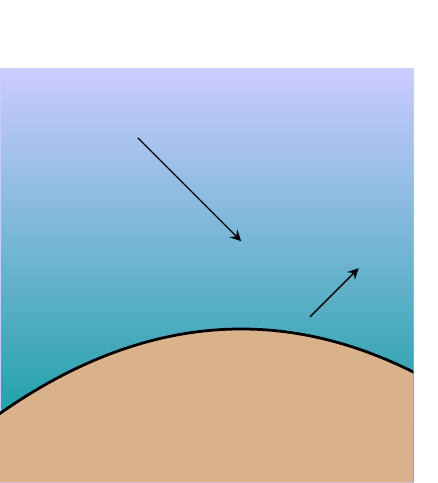}
\end{center}
\end{figure}

The proposed inversion strategy has two main components. The first is the parametrization of the seafloor domain and the corresponding physical models of the wave propagation and backscatter, described in \S \ref{sec:forward}. Simulations corresponding to a seafloor parameter space are used to generate a library of acoustic templates, described in \S \ref{sec:library}. The second component, discussed in \S\ref{sec:inverseproblem},  is the optimization procedure for matching predictions of the seafloor parameters with given backscatter data. Numerical simulations are presented in \S \ref{sec:numerics}, and we conclude in \S \ref{sec:conclusions}.

\section{Forward model}\label{sec:forward}

In high frequency wave propagation in an underwater acoustic environment, the acoustic penetration of the seafloor is limited and it suffices to only consider the scattering effects of the upper sedimentary layer (see Figure \ref{fig:SFdiagram}). Here, we focus on the two-dimensional setting though these approaches can be generalized to three spatial dimensions. The horizontal coordinate is given by $x$  and the vertical coordinate is given by $y$.

\begin{figure}\caption{ In the domain $\Omega$, sources and receivers are located near the sea surface generate an acoustic signal that acts as an incident plane wave with incident angle $\alpha$ in the region $D$ containing the water-sediment interface. Interactions with the sediment layers cause scattered waves in the direction $\beta$.} \label{fig:fwd}

\begin{center}
\includegraphics[scale=.8]{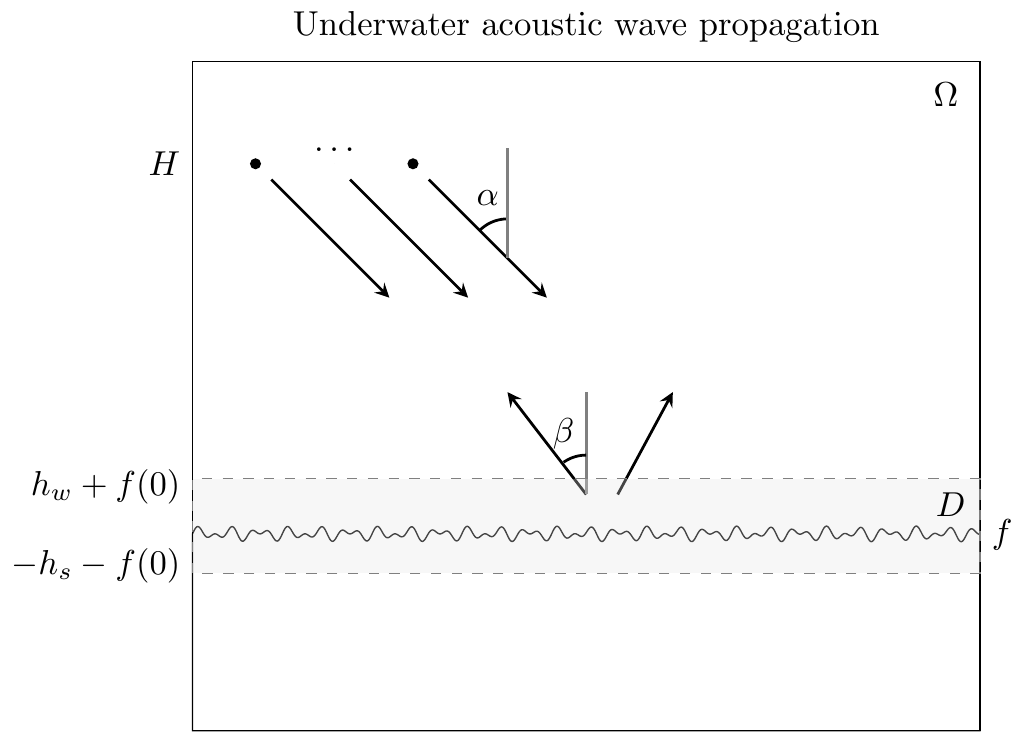}
\end{center}
\end{figure}

We will consider the problem of recovering the fine scale geometry and material type of the seafloor located at $y=f(x)$, assuming large scale variations are known. The scattering process is modeled in the rectangular domain $\Omega$, and the $n$ measurement locations at a depth $y=H$ near the sea surface are given by the set $X$, \begin{align} X= \{(x'_{j},H)\mid 1\leq j \leq n\} \labeleq{X}.\end{align}

\noindent Because the most significant contributions to the measured signal are produced from scattering near the ocean bottom, high fidelity modeling is needed in a thin layer $D\subset \Omega$ containing the water-sediment interface,
\begin{align*}
D=\Omega\cap\{-h_{s}-f(0)\leq y \leq h_{w}+f(0)\}.
\end{align*}
Model considerations are used to determine the vertical thickness of the truncated sediment layer and water layer, given by $h_{s}$ and $h_{w}$, respectively. For example, in our simulations the value of $h_{s}$ was chosen to account for the attenuation of the acoustic energy in the sediment layer and prevent artificial reflections. The value of $h_{w}$ was chosen so that oscillations from the seafloor could be accurately modeled in the seawater domain using the technique of numerical microlocal analysis. Figure \ref{fig:fwd} shows a diagram $\Omega$ with the thin layer $D$ highlighted in yellow. 

\subsection{Governing equation}

As in the majority of fluid-sediment modeling \cite{Chapman2001, Jackson2007,  Jackson1986}, we treat the sediment layer as a fluid and model time-harmonic acoustic pressure waves as solutions $p$ to Helmholtz equations with varying coefficients,
\begin{align}
{\rho_{0}}\nabla\cdot\(\frac{1}{\rho_{0}} \nabla p\)  +\frac{\omega^{2}}{c_{0}^{2}} p &=0 & \text{ in }  D\cap \{ y>f(x)\}, \labeleq{fwd1}\\
p \vert_{+} = p\vert_{-}, \qquad \frac{1}{\rho_{0}}\frac{\partial p}{\partial \vec{n}}\vert_{+}& =\frac{1}{\rho}\frac{\partial p}{\partial \vec{n}}\vert_{-}  &\text{ on }  D\cap \{ y=f(x)\},\labeleq{fwd2}\\
{\rho}\nabla\cdot\(\frac{1}{\rho} \nabla p\)  +\frac{\omega^{2}}{c^{2}} p +i\gamma {\omega}{p} &=0  & \text{ in }  D\cap \{ y<f(x)\},\labeleq{fwd3}
\end{align}
where $\omega$ is the frequency variable. The model coefficients represent the sound speed $c$ and material density $\rho$, given in the water layer by $c_{0}=1030$ m/s$^{2}$ and $\rho_{0}=1500 $kg/{m$^{3}$. The attenuation of high frequency acoustic waves in the sediment is given by $\gamma$, and we assume that the attenuation in water is negligible. 

There are two coupling conditions on the interface in \refeq{fwd2}. The first is a continuity of pressure along the fluid-sediment interface. The second condition is the continuity of displacement normal to the boundary of the interface. Here the vector $\vec{n}$ is normal to the curve given by $(x, f(x))$ and $\vert_{+}$ and $\vert_{-}$ refer to the limits of a function or derivative as $y$ approaches $y=f(x)$ from above and below, respectively. As a remark, in situations where it is desirable to model more that one sedimentary layer, this model can be augmented with additional interface conditions of the type \refeq{fwd2} and equations of the type \refeq{fwd3}.

For boundary conditions, the fast attenuation in the sediment layer is modeled with the condition $p=0$ on the lower boundary of the domain, and absorbing boundary conditions representing the incoming wave are imposed on the upper boundary. In general, the conditions on the left and right boundaries are determined by the assumptions on scattering contributions from neighboring sediment layers. 

\subsection{Problem setup}
The process begins when an acoustic pulse, called a ping, is generated by a source near the surface. Close to the seafloor, the signal acts as an incident plane wave traveling in the direction given by the vector $k_{\alpha}=(k \cos \alpha, -k \sin \alpha)$, where $\alpha$ is the incident angle and $k=\frac{\omega}{c}$ is the wavenumber (see Figure \ref{fig:fwd}). Supposing that the strength of the source is $P_{0}$, the incident wave $p_{inc}$ in the domain $D$ is expressed as  \begin{align}
p_{inc}(x,y, \alpha) = P_{0}e^{ i  k_{\alpha} \cdot (x,y)}.\labeleq{pinc}
\end{align}

The interaction of these waves with the seafloor produces outgoing waves in multiple directions depending on the roughness and material type of the sediment layer. Sufficiently far from caustics, the scattered wavefield $p_{out}$ can be approximated by the superposition given by
\begin{align}
p_{out}(x,y)\simeq \sum_{\beta}P_{\beta}(x,y) e^{ i k_{\beta}\cdot (x,y)},\labeleq{wavedec}
\end{align}
where $P_{\beta}(x,y)$ is the complex amplitude of the outgoing wave in the direction $\beta$. 

Acoustic backscatter is measured in most sonar applications in which the receivers and source locations are the same. The receivers record the strength of the outgoing wave traveling toward the source in the direction $\beta=\alpha$. Here we are not focusing on technologies that involve distance measurements, and model the received signal $Y\in \RR^{n}$ as measurements of the backscatter strength on the domain $X$ in \refeq{X}, \begin{align}
Y_{j}=|P_{\alpha}(x'_{j}, H)| \qquad 1\leq j\leq n \labeleq{backscattered}. 
\end{align}
\noindent  The magnitude of the coefficient $P_{\alpha}$ gives the strength of the backscattered wave in the direction $\alpha$. 

The typical spatial domain $\Omega$ is large, on the order of kilometers, while the acoustic wavelength is a few centimeters. A domain of size 100 m covers thousands of times the acoustic wavelength of a 20 kHz source signal in range, and tens of wavelengths in depth \cite{Fox1985}. Full resolution of the signal in the water layer alone would require a grid of at least five points per wavelength, on the order of hundreds of millions of points for 2D resolution. Several predictions of the forward model must be made in the template matching process, and the computational task quickly becomes too expensive. 

The original model can be reduced by dividing the domain $D$ into $N$ disjoint segments of length $\Delta s$, and determining the backscattering caused by the seafloor occupying the $i^{th}$ segment of the seafloor domain, 
\begin{align}
D^{i}=D\cap \{(i-1)\Delta s<x< i\Delta s\},\quad 1\leq i \leq N \labeleq{Di}. 
\end{align}
The width $\Delta s$ of each segment is chosen to be roughly $10-20$ times the acoustic wavelength. Then, we can reasonably assume that the scattering features of interest in each segment of length $\Delta s$ is independent from segments that are spatially far away and approximate the backscattering using only local simulations. We will then design a procedure for approximating the backscattering $\refeq{backscattered}$ on small domains, each of size $\Delta s$ using a 1) parametrization of the seafloor and 2) localized simulations.

\subsection{Seafloor parametrization} \label{sec:parameterization}
A simplified description of the seafloor domain can be made by assigning a vector of parameters $m^{i}$ that describe the seafloor properties within segment $D^{i}$. The identification of seafloor parameters that can be recovered from backscattered data can enable efficient predictions of the scattering effects from a wide range of sea floors. Seafloor parameters fall into three main categories: experimental parameters, geometric parameters, and geoacoustic parameters. Experimental parameters $m_{E}$ are determined from the physical scenario, and it will be assumed that these parameters are known beforehand. Within the segment $D^i$ it can be assumed, for our purposes, that the geometric parameters $m_{G}$ can be modeled as smoothly varying functions that take a continuous range of values. Though the full seafloor domain consists of many sediment types, we assume that each seafloor patch, when it is small enough, consists of a uniform sediment, that is,  $m_{A}$ remains constant. The parameter $m_{A}$ is qualitative; we will assign to the segment $D^i$ one label from a list of possible material types. Descriptions of the categories are below:

\begin{itemize}

\item {\bf Experimental parameters, $m_{E}$.} These parameters include the locations of the source and receiver,  source strength, source frequency, and  incident angle. The source and receiver locations are colocated near the surface of the water, and the strength of the acoustic source is a constant and given by the sonar manufacturer. For sidescan sonars this is typically between $200-230$ decibels \cite{Pailhas2010}. High frequency signals produce strong scattering from artifacts in the seabed that have variations on the scale of the acoustic wavelength. These artifacts cannot be easily detected using lower frequencies. We consider source signal s with a frequency of 20 kHz. The acoustic wavelength is then on the scale of centimeters. The range of the incident angle $\alpha$ is determined by the instrument. Specifically, in sidescan sonars, $\alpha$ lies within the critical range of $30^{\circ}\pm 15^{\circ}$.

\item {\bf Geometric parameters, $m_{G}$}. The interface function $f(x)$ depends on fine scale features, such as the surface roughness characteristic, and large scale features, such as the overall curvature of the interface. Statistical models of seafloor roughness parameterize the distribution of the normalized amplitude of representative component sinusoids \cite{Fox1985}. The idea is to mimic a realistic seabed that is not necessarily periodic. The parametrization of $f$ can be more general, but in this study, we model a seabed with sinusoidal ripples as an example,
\begin{align}
m_{G}\rightarrow f(x, m_{G}) =  .01\(\sin(2\pi m_{G1}x)+m_{G2}\sin(2\pi m_{G3}x)\) \labeleq{f}.
\end{align} 
The geometric parameter $m_{G}=(m_{G1},m_{G2}, m_{G3})$ refers to the frequencies $m_{G1}$ and $m_{G3}$ and aspect ratio $m_{G2}$ of the component sinusoids.
The choice of spatial frequencies $m_{G1}\simeq 15$ and $m_{G3}\simeq 25$ produces a spatial variations on a scale similar to the characteristic wavelength of the incoming acoustic wave (see Figure \ref{fig:fwd}). We chose values for the aspect ratio that lie in the interval $.5\leq m_{G2} \leq 1$.

\item {\bf Geoacoustic parameters, $m_{A}$.} The acoustic description of the the sound speed, density, and attenuation, is purely based on material type of the sediment layer. The influence of the material type is not explicitly used in the interface function, however it is in practice reflected in the fine scale ranges of the acoustic response. For example, the fine scales for mud, sand, sea grass, and rock are of an increasing order in their characteristic length scales.

Here, we included sediment types with distinctive geoacoustics. The sound speed in fine sediments such as clay are typically ``slow'' (1-5\% less than water), while the speed in sand is ``fast'' (10-15\% more than water) \cite{Jackson2007}. In all of these sediments, the sound attenuation $\gamma$ is roughly 0.5 dB per 1 meter per 1 kHz. We chose not to distinguish attenuation based on material type.  In Table $\ref{table:mat}$,  the density $\rho_{m_{A}}$ and sound speed $ c_{m_{A}}$ of the different materials used in this study is given. 
\begin{table}
\caption{Geoacoustic parameters}\label{table:mat}
\begin{center}
\begin{tabular}{|c|c|c|}
\hline
Material type & density (kg/m$^{3}$)& sound speed (m/s)\\
\hline
sand& $\rho_{sand}=2000$& $c_{sand}=1668$\\
clay& $\rho_{clay}=1170$& $c_{clay}=1518.9$\\
rock & $\rho_{rock}=2870$& $c_{rock}=6000$\\
metal & $\rho_{metal}=8050$ & $c_{metal}=6100$\\
\hline
 \end{tabular}
\end{center}
\label{default}
\end{table}%

%

The possible material types $m_{A}$ that represent the sediment layer are 
\begin{align*}
m_{A}=\text{`sand', `clay', `rock', or `metal'}.
\end{align*}
A typical sediment layer is depicted in Figure \ref{fig:materialtype} (a). The material type $m_{A}=$ `metal' is used to label a sand layer containing a buried rectangular metal object, however this can easily be extended to objects buried in other sediments or objects that are only partially buried. The metal domain, depicted in Figure \ref{fig:materialtype} (b), is $D' \cap \{ y < -\epsilon_{f}\}$. The object lies purely beneath the interface, $\text{inf}_{x}f(x)+\epsilon_{f}>0$. The depth parameter $\epsilon_{f}>0$ is chosen to be large enough so that the object can be detected using high frequency waves. If $\epsilon_{f}$ is too large, other techniques must be used, for example, using lower frequency pings. In this study we only considered objects buried in a sediment layer of sand, however this can easily be extended to objects buried in other sediments or objects that are only partially buried.
 
In summary, the sediment layer $D'\cap \{ y <f\}$ is assigned the geoacoustic properties $m_{A}\rightarrow (\rho_{m_{A}}, c_{m_{A}})$ for $m_{A}=\text{sand, clay, rock,}$. For $m_{A}=\text{metal}$, these properties are given by
\begin{align}
m_{A}\rightarrow\begin{cases}(\rho_{sand}, c_{sand})& -\epsilon_{f}< y ,\\
(\rho_{metal}, c_{metal})&  y < -\epsilon_{f}
\end{cases}.\labeleq{mA}
\end{align}

 \begin{figure}\caption{Diagrams depicting typical sediment layers modeled using acoustic templates. Shown trom left to right: a homogeneous sand layer (a), a sediment layer containing a buried metal object (b), and a sediment layer containing a horizontal transition between two materials (c).  }\label{fig:materialtype}
\begin{center}

\includegraphics[scale=1]{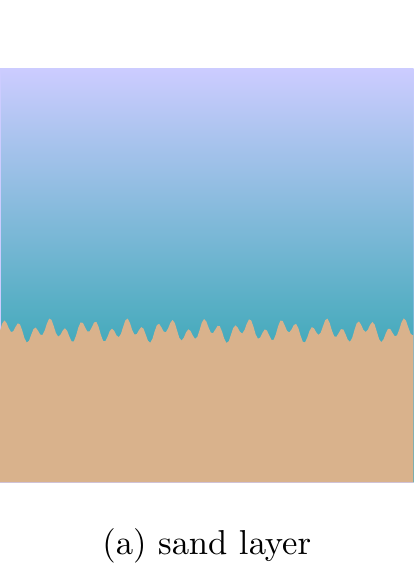}\hspace{1em}\includegraphics[scale=1]{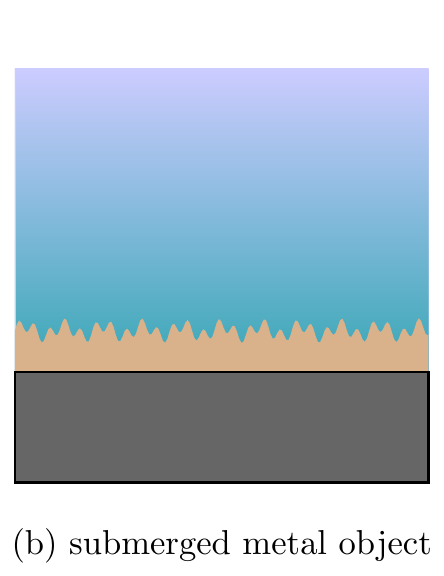}\hspace{1em}\includegraphics[scale=1]{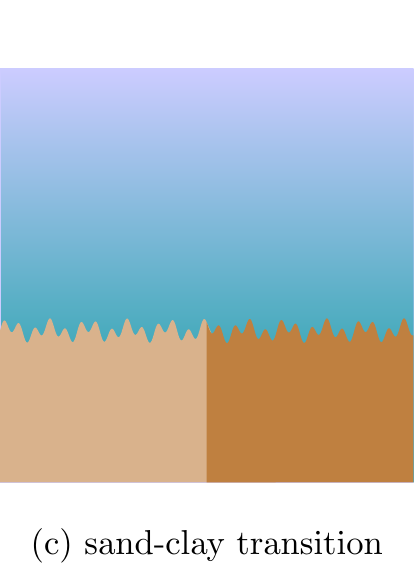}
\end{center}
\end{figure}

\end{itemize}

We chose a parameter space $\mathcal{M}$, consisting of acoustic and geometric parameters,
\begin{align*}
\mathcal{M}=\{m=(m_{A}, m_{G})\},
\end{align*}
where $m_{A}$ is a parameter that determines the sound speed and density of the sediment and $m_{G}$ determines the interface function.

\subsection{Localized modeling}\label{sec:localized}
For large seafloor domains, simulations of high frequency Helmholtz equations are often cost-prohibitive. A key component of our framework involves determining the acoustic response of smaller spatial segments of seafloor domain by matching the observed acoustic response with a collection of precomputed signals. These signals, called acoustic templates, represent the backscattering from small patches of the sediment layer with the typical geometric and material characteristics. Based on the previous discussion, a vector of parameters $m=(m_{A}, m_{G})$ can be used to describe each segment $D^{i}\subset D$. Forward simulations can be used to determine the contribution to the backscatter from the parameterized seafloor.

Let $D'\subset D$ be a generic template domain of width $\Delta s$, and for $m\in \mathcal{M}$, define the interface function $f$ determined by $m_{G}$ and the formula \refeq{f}. In the sediment layer $D'\cap\{y<f\}$, the density and sound speed of the material are determined by $m_{A}$ as in \refeq{mA}. Then, the solution $\hat{p}_{m}$ to the Helmholtz equations $\refeq{fwd1}-\refeq{fwd3}$ approximates the scattering from the seabed. Let $\hat{b}_{m}$ be the backscatter computed from $\hat{p}_{m}$ in the domain $D$ given by \refeq{backscattered}. Then, the acoustic response for each segment of the seabed has the correspondence
\begin{align}
m\rightarrow \hat{b}_{m}, \qquad (\hat{b}_{m})_{j}=|P_{\alpha}(x'_{j}, H)|, \quad 0\leq x'_{j}<\Delta s \labeleq{ym}.
\end{align}

The boundary conditions in the local model $\refeq{fwd1}-\refeq{fwd3}$ must be chosen carefully. As a result of our scale assumptions, the vast majority of segments $D^{i}\subset D$ have the same material type as the neighboring segments $D^{i-1}$ and $D^{i+1}$. For simplicity, we impose periodic boundary conditions that account for the scattering contributions from neighboring domains of the same material type and geometry. This is a common strategy in scattering theory in the context of micro-optics, where periodic structures are called diffraction gratings \cite{Bao1995, Bao2005, Bao2013, Bruno1991}.

On the other hand, there are instances where segments of the domain differ in material type from neighboring segments. For each $i=1,\hdots N-1$, we will call the pair of material types in the cells $D^{i}$ and $D^{i+1}$ a ``transition'' at the location $x=i\Delta s$. For example, the sediment layer shown in Figure \ref{fig:materialtype} (c) contains a `sand-clay' transition at location $x=i \Delta s$. Sediment layers containing a transition can be modeled by embedding the template domain into a larger domain in which periodic boundary conditions are imposed.

The sources of modeling error include the local representation of the backscatter using truncated domains, the smoothness assumptions on the parameters, and any resonance phenomena occurring from the choice of periodic boundary conditions.

\subsection{Numerical microlocal analysis for solutions to the Helmholtz equation}\label{sec:mla}
Due to the fine scale structure of the seafloor, interactions between the sediment layer and the the incoming acoustic wave produce complex scattering effects. The modeling of synthetic backscatter requires accurate details about the propagation of outgoing waves, especially the waves that propagate along the receiver direction. This motivates the need for strategies for extracting information from the solution $p$ to the Helmholtz equation $\refeq{fwd1}-\refeq{fwd3}$ that represents the backscattered signal  $Y\in \RR^{n}$ given by \refeq{backscattered}. 

Common approaches for modeling the scattering effects are based on geometric ray tracing or full solutions to wave equations \cite{Bell1997, Elston2004}. Approaches such as Kirchhoff approximation (requiring the seafloor interface to be smooth) and small perturbation theory (requiring the roughness scale to be much smaller than the acoustic wavelength)  are independent of the wave theory used and do not account for volume scattering within the sediment \cite{Groen2010, Jackson1986}. Furthermore, many models for seafloor backscatter are often developed for only seabeds consisting of isotropic roughness.

Here, we use a different tool that is sensitive to the underlying seafloor material properties as well as microtopography.  Our proposed strategy has two components. The first component involves modeling scattering by interpreting solutions to Helmholtz equations as superpositions of waves traveling in different directions. The second component involves modeling the backscattered signal  $Y$ and extracting information about the waves traveling in the backscattered direction.

\begin{figure} \caption{ Plot of the Helmholtz solution to $\refeq{fwd1}-\refeq{fwd3}$ with incident angle $\alpha=\pi/6$ and varying material types. Indicated in black are the circles with radius $r_{0}$ centered at the locations $x=.5, x=3, $ and $x=6$. }\label{fig:microlocalsoln}
\vspace{2em}
\includegraphics[scale=1]{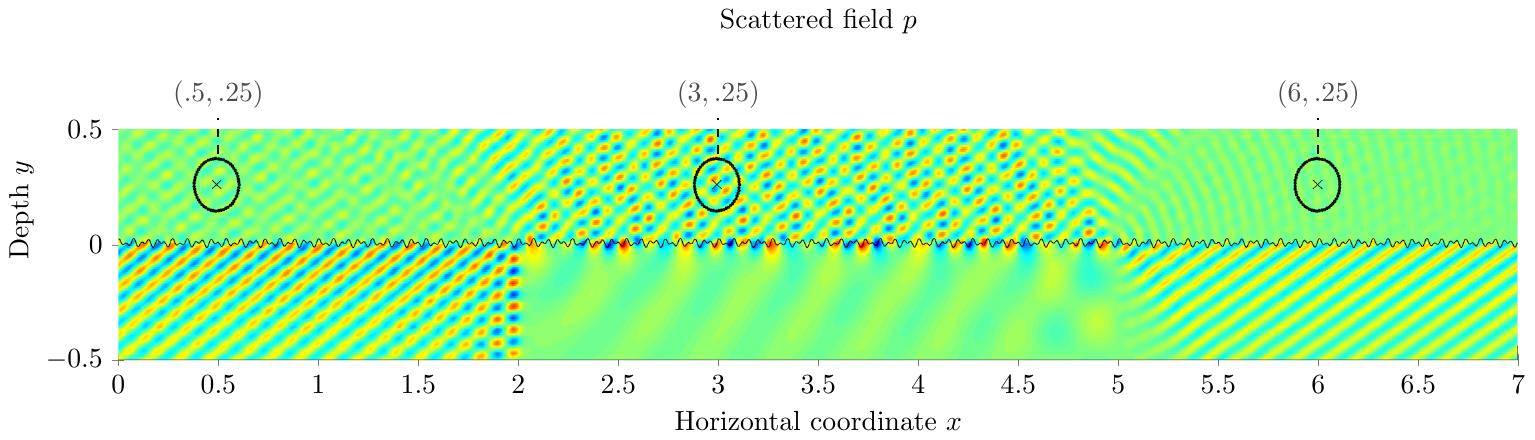}
\end{figure}

The measured scattering from the seabed is approximated using numerical microlocal analysis of solutions to Helmholtz equations, a technique for determining the number, directions, and the corresponding complex amplitudes of the crossing rays in domains where geometric optics is valid. Described in detail in \cite{Benamou2004, Landa2011}, the main idea is that if observation points $(x_{0},y_0)$ are far away from caustics and focus points and the wavelength $\lambda_{0}=2\pi/k$ contains variations on a scale sufficiently small with respect to the geometry, the solution $p(x,y)$ to $\refeq{fwd1}-\refeq{fwd3}$ behaves like a superposition of plane waves. Adopting the notation in \cite{Benamou2004, Landa2011}, the local plane wave approximation of $p(x,y)$ restricted to a circle centered at $(x_{0},y_0)$ with sufficiently small radius $r_{0}$ is
\begin{align}
p(x,y)\simeq \sum_{l=1}^L P_{\beta_{l}}(x,y)e^{ik\phi_{l}\cdot (x,y)},\labeleq{pwansats}
\end{align}
when $|(x,y)-(x_{0},y_0)|<r_{0}$ (Figure \ref{fig:microlocalsoln}). The values of $ P_{\beta_{l}}(x,y)$ correspond to the complex amplitudes of the outgoing waves traveling in the directions  $\phi_{l}=(\cos(\beta_{l}), \sin(\beta_{l}))$. The function $P_{\beta_{l}}(x, y)$ is locally smooth and independent of wavelength.

The algorithm recovers the number $L$ of crossing rays and the magnitudes $|P_{\beta_{l}}(x_{0},y_0)|$ of the scattered waves in the direction given by $\beta_{l}$. In the derivation, the 2-D Jacobi-Anger expansion is used to represent the local solution as a series involving Bessel functions, which is then appropriately truncated at a threshold  $\hat{L}(r_{0})=r_{0}+5r_{0}^{1/3}$. Finally, the desired $L=2\hat{L}(r_{0})+1$ unknown amplitudes are determined using a Fourier-type inversion formula involving the Bessel function of order $l$, denoted by $J_{l}(r_{0})$.  Regularization can be applied to improve the stability of the reconstruction formula.

The following steps summarize the algorithm for approximating the backscattering at each observation point $(x_{0},y_0)$:

\begin{description}
\item{\bf Step 1.} Define $L$ equidistant points on the interval $[0, 2\pi)$, 
\begin{align*}
\beta_{l}=\frac{2\pi({l}-1)}{L}, \qquad {l}=1, \hdots, L,
\end{align*}
and let $\{\tilde{P}_{l}\}$ be the grid function of sampled values of the Helmholtz solution $p(x,y)$ on the restricted circle centered at the observation point $(x_{0},y_0)$ with radius $r_{0}/k$, $r_{0}=3\pi$ and angles $\beta_{l}$, 
\begin{figure} \caption{ The polar plots indicate the direction  and intensity of the scattered waves at the observation points indicated in Figure \ref{fig:microlocalsoln}. }\label{fig:mlapolar}
\begin{center}

 \includegraphics[scale=.8]{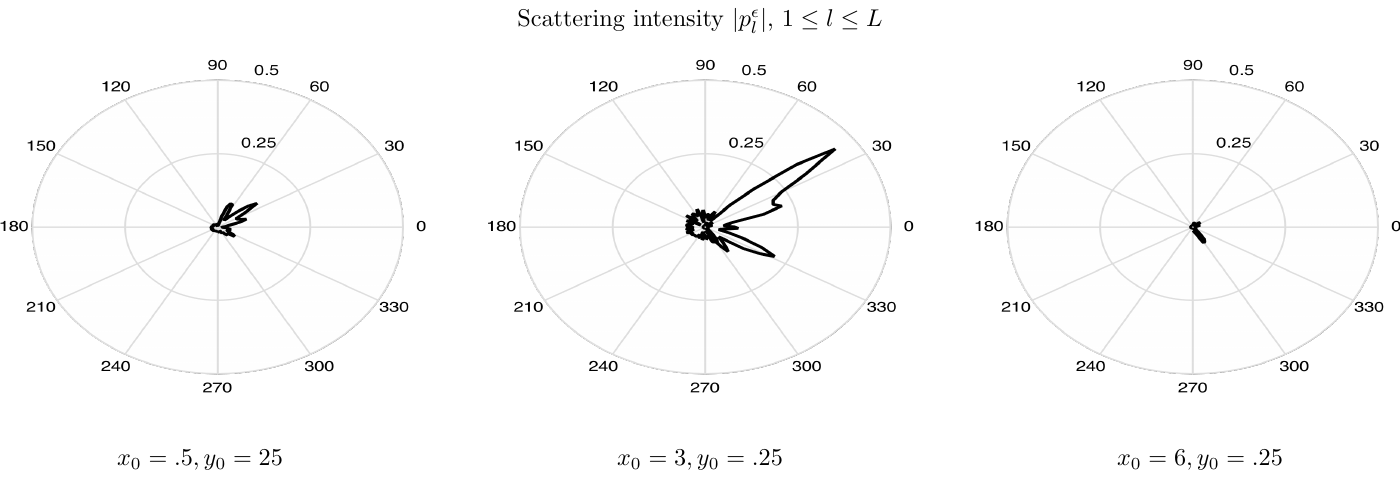}

\end{center}
\end{figure}

\begin{align*}
{\tilde{P}_{l}}=p\( x_{0}+\frac{r_{0}}{k}\cos(\beta_{l}) ,y_0+\frac{r_{0}}{k}\sin(\beta_{l})\),\qquad {l}=1, \hdots, L.
\end{align*}
\item{\bf Step 2.} Let $\epsilon>0$ be a regularization parameter, and denote by $FFT$ and $FFT^{-1}$ the discrete Fourier transform and the inverse discrete Fourier transform, respectively. Define 
\begin{align*}
\{p_{l}^{\epsilon}\}=2\pi FFT^{-1}\left\{\frac{FFT\{\tilde{P}_{l}(2\hat{L}(r_{0})+1)J_{l}(r_{0})\}}{i^{q}[(2\hat{L}(r_{0})+1)^{2}J_{q}(r_{0})^{2}+4\epsilon\pi^{2}]}\right\}
\end{align*}
\item{\bf Step 3.} $|P_{\beta_{l}}(x_{0},y_0)|\simeq|p_{l}^{\epsilon}|$.

\end{description}

Figure \ref{fig:mlapolar} shows the resulting polar plots of $|p_{l}^{\epsilon}|$ taken at the three observation points indicated in Figure \ref{fig:microlocalsoln}. Figure \ref{fig:mlaAll} gives a surface plot of the scattering $|p_{l}^{\epsilon}|$ at measurements taken across the full spatial domain (horizontal axis) and full range of scattering directions $\beta$ (vertical axis).

In the sonar applications of present interest, receivers only detect the strength of signals propagating in the backscattered direction $\beta=\alpha$. We estimate this quantity using linear interpolation between values in the given pairs $(\beta_{l}, p^{\epsilon}_{l})$. For $\beta$ is in the interval $[\beta_{l}, \beta_{l+1})$, the amplitude of the backscattered wave \refeq{backscattered} is approximated by

\begin{align*}
|P_{\beta}(x_{0},y_0)|\simeq p_{l}^{\epsilon}+(p_{l+1}^{\epsilon}-p_{l}^{\epsilon})\frac{\beta - \beta_{l}}{\beta_{l+1}-\beta_{l}}.
\end{align*}

\begin{figure}\caption{Micro local analysis is applied to solutions to $\refeq{fwd1}-\refeq{fwd3}$ and a plane wave source with incident $\alpha=\pi/6$. The scattering strength  in each direction is calculated using numerical microlocal analysis. }\label{fig:mlaAll}

\begin{center}
\includegraphics[scale=1]{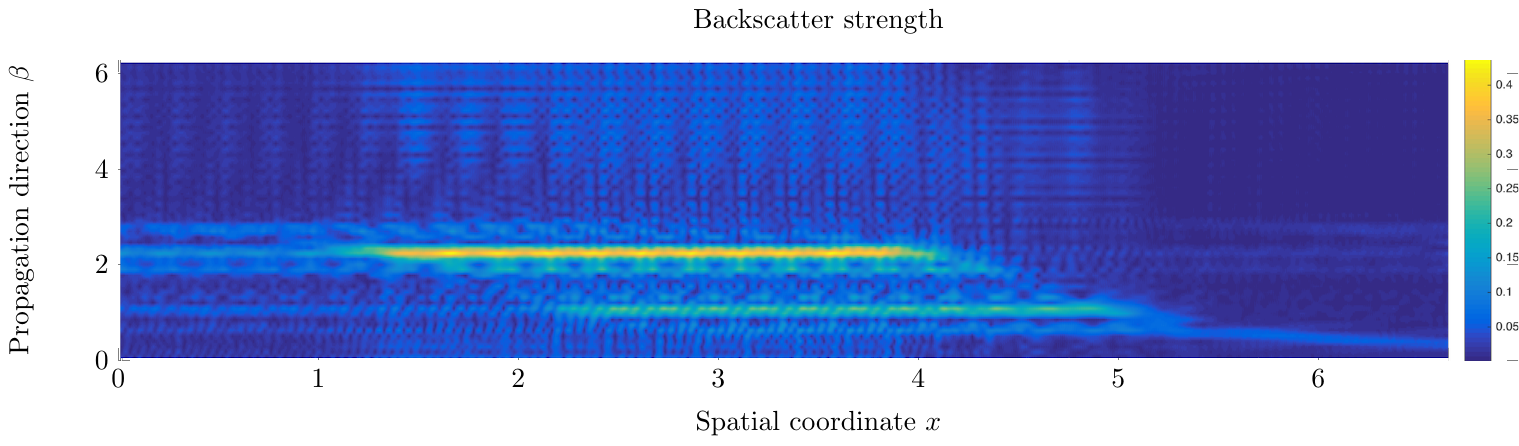}

\end{center}
\end{figure}

\section{Construction of the library of acoustic templates}\label{sec:library}

%
%

The library $\mathcal{L}$ of acoustic templates, is generated from solutions of $\refeq{fwd1}-\refeq{fwd3}$, 
\begin{align}
\mathcal{L}=\{\hat{b}_{m} \mid m\in \mathcal{M}\}, \labeleq{library}
\end{align}
where $\hat{b}_{m}$ is a vector of the backscatter \refeq{backscattered} measured at locations $(x_{j}', h)$ in the template domain  $D'$.  The key to successful inversion lies on the customization of the library $\mathcal{L}$, especially with respect to surface geometry and material properties; combinations of  templates from a well-designed library can be used to generate acoustic signals for a range of scenarios. 

Careful thought must be put into determining the resolution of the library. The template-matching technique is largely dependent on the number of signals that must be generated in the library; the processing load increases with the number of templates required, and computations may become cost-prohibitive. However, in some cases, small changes in seafloor parameters result in highly variable scattering patterns, and if the sampling rate is too low, critical features may not be resolved. Justifying our choices using a constraints of the sonar imaging problem, numerical experiments, and heuristics, we set a reasonable range and resolution for parameters included in the library. 

The four sediment types were chosen to reflect materials with distinctive geoacoustic properties. We also added redundancy to the library by including multiple templates for each material type with consideration to the neighboring sediments, with the goal of object detection and improved classification. The templates were generated by solutions to $\refeq{fwd1}-\refeq{fwd3}$, first in the template domain with periodic boundary conditions, and then in a larger domain containing a transition of different sediment types. With such an enriched library comes enhanced recovery of material type, lowering the rate of false object detection. Numerical simulations and details will be given later.

The range of the geometric parameters is chosen fine enough to demonstrate smooth variations in the signal. Figure \ref{fig:mg-sensitivity} displays the sensitivity of the backscattered data with respect to the parameters $\alpha$, $m_{G1}$, $m_{G2}$ and $m_{G3}$. The sediment layer consists of sand in all four cases. We decided a resolution of $N=20$ is needed to capture distinct characteristics in the spatial frequencies.   Plots suggest $m_{G2}$ has some regularity and can be predicted by neighboring values, therefore, we determined that a resolution of $N=4$ is needed.  These results are summarized in Table \ref{table:mg-sensitivity}.

\begin{figure} \caption{The surface plots show variations in the backscatter $\hat{b}_{m}$ across the spatial domain (horizontal axis) and the parameter domain (vertical axis).}
\begin{center}
 \resizebox {.4\columnwidth} {!} {
\includegraphics[scale=1]{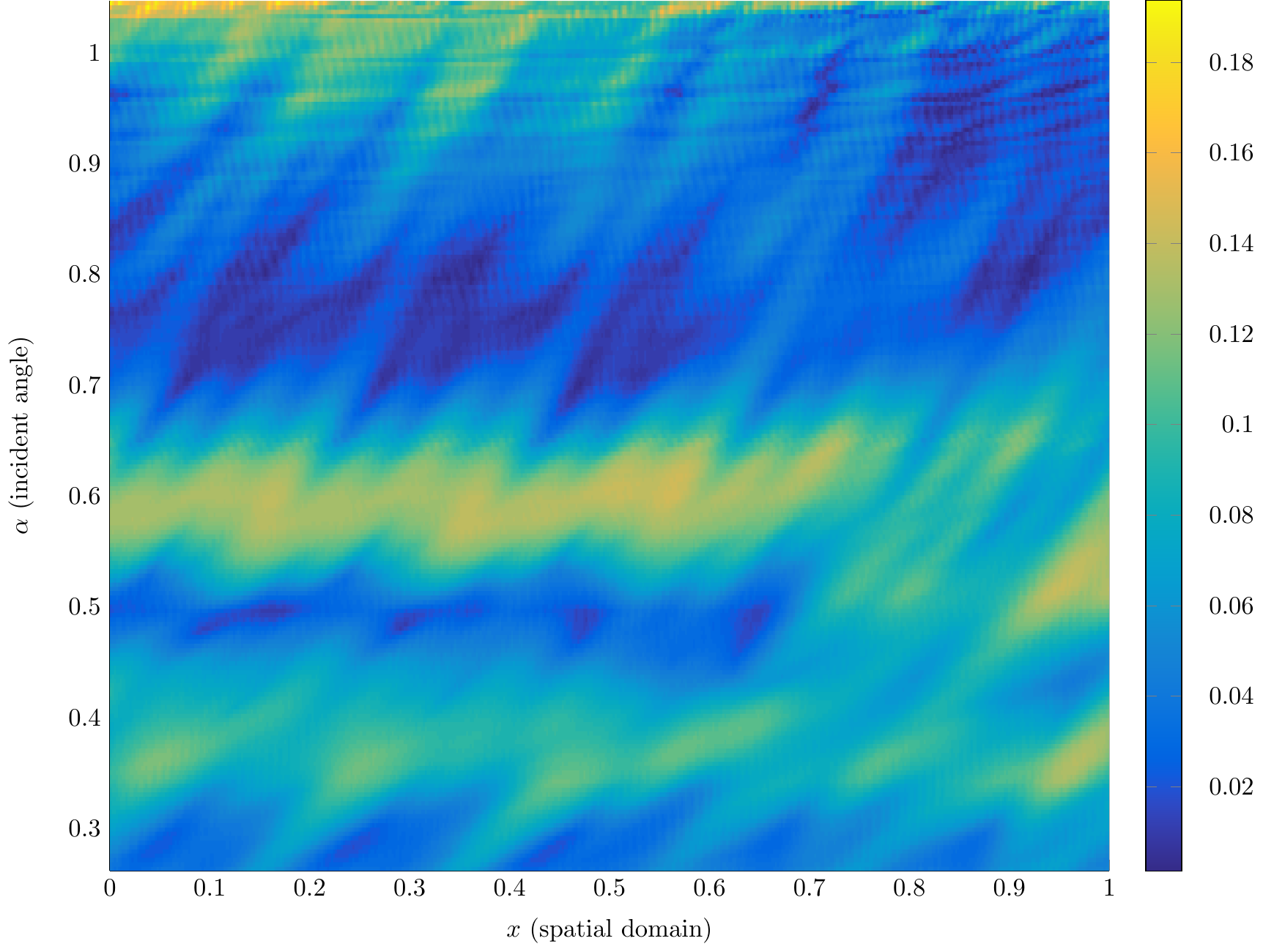}}
 \resizebox {.4\columnwidth} {!} {
\includegraphics[scale=1]{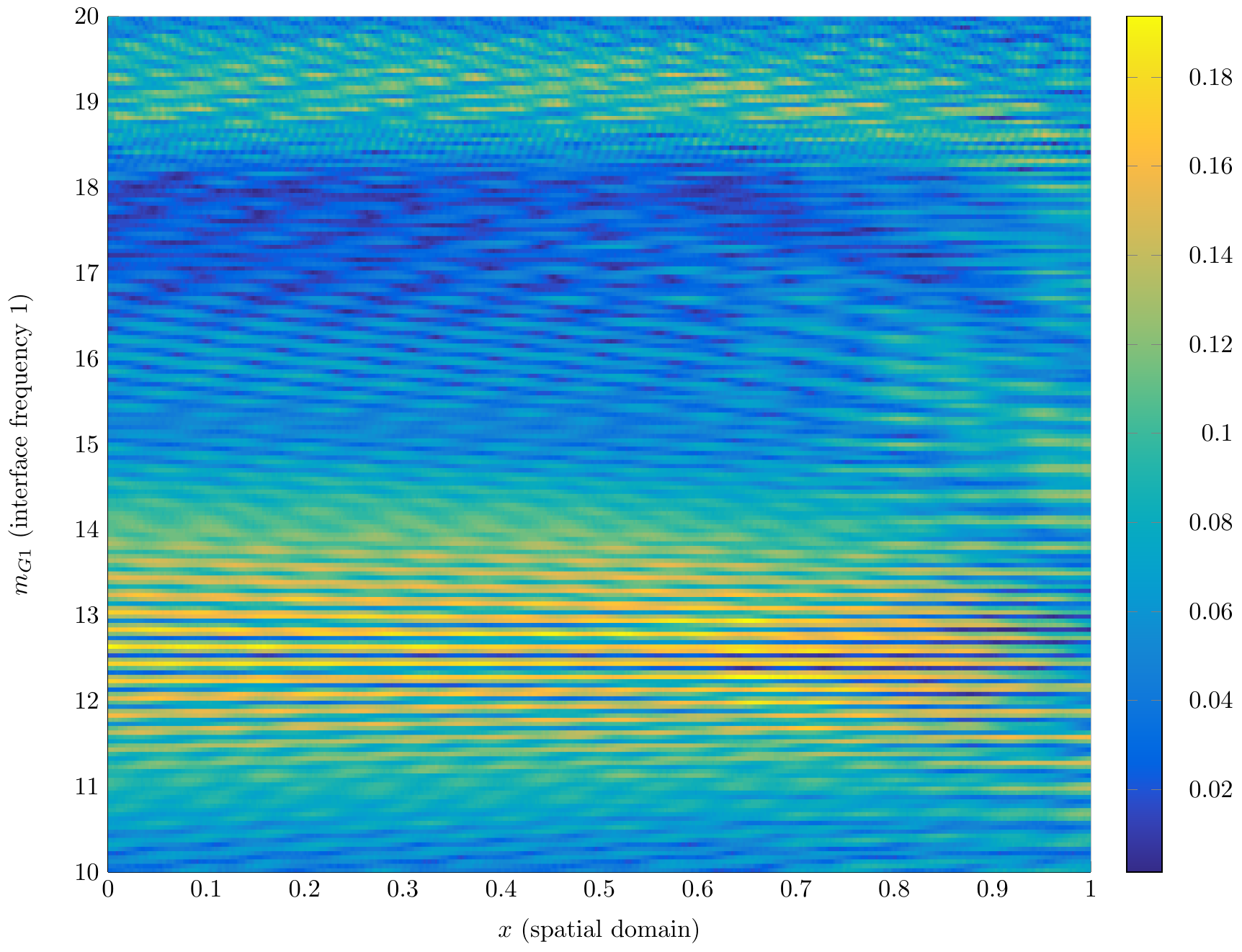}}
 \resizebox {.4\columnwidth} {!} {
\includegraphics[scale=1]{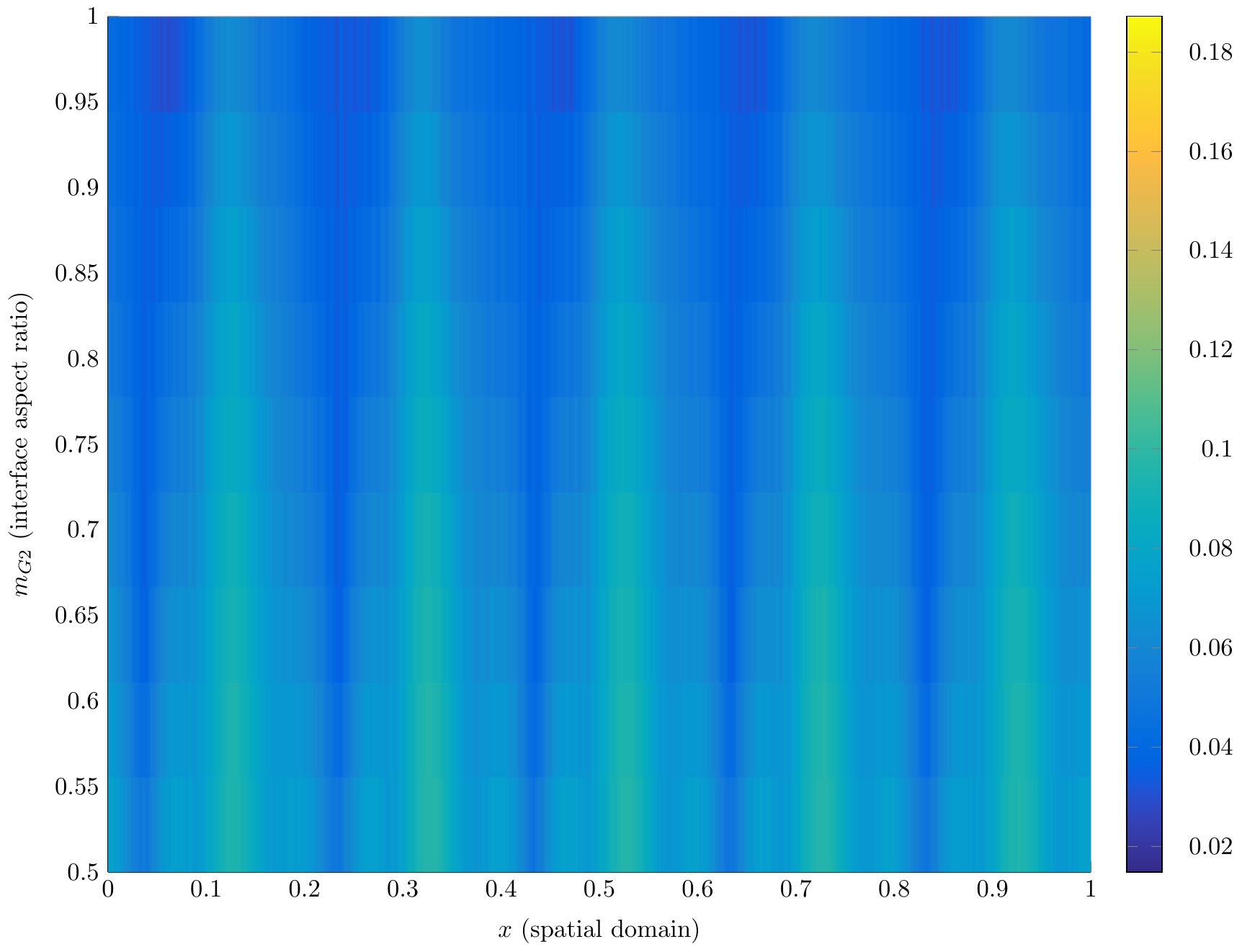}}
 \resizebox {.4\columnwidth} {!} {
\includegraphics[scale=1]{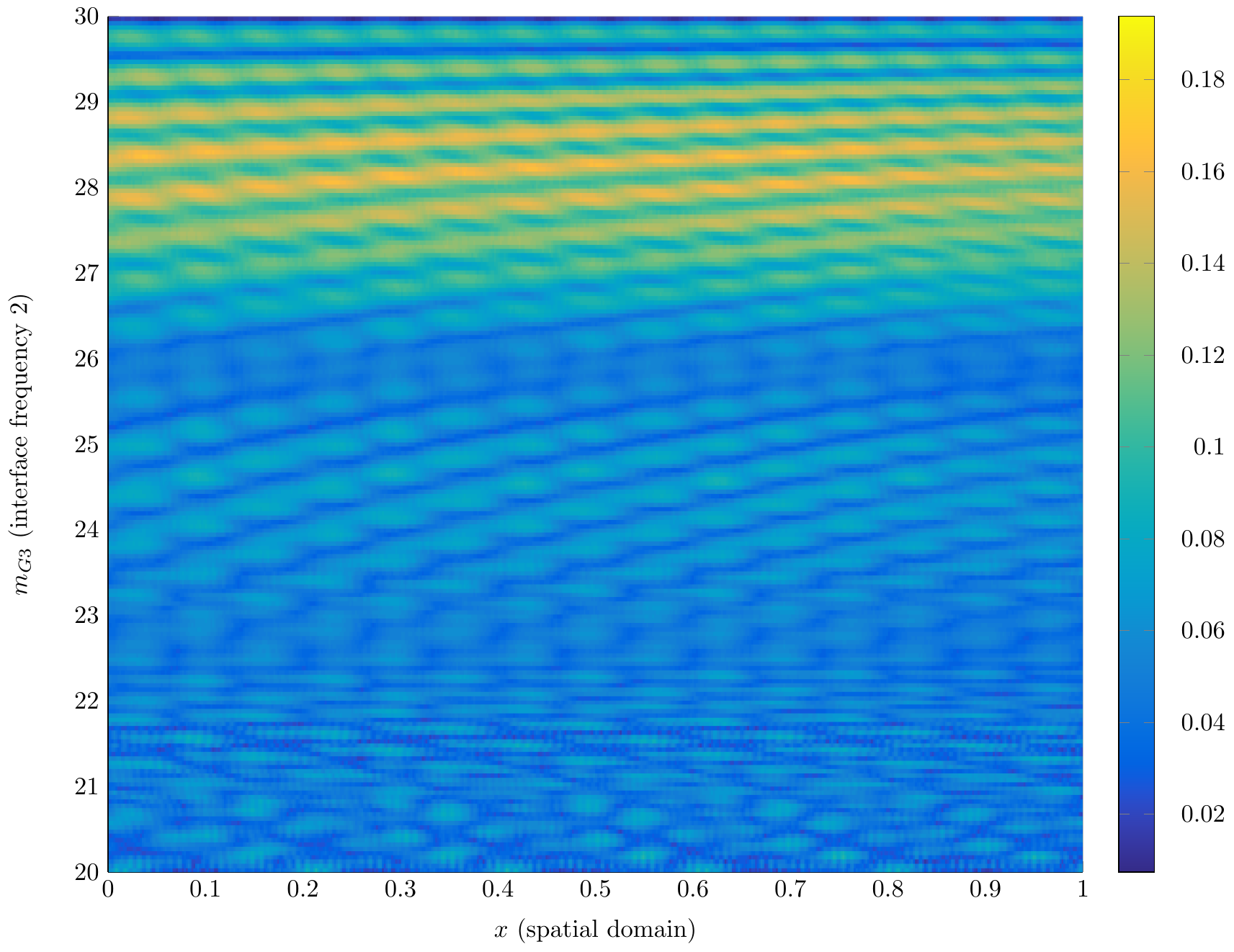}}
\end{center}
\label{fig:mg-sensitivity}

\end{figure}

\begin{table}
\caption{The table below gives the range and resolution of geometric parameters $m_{G}$ chosen for the chosen library of templates.}\label{table:mg-sensitivity}
\begin{center}
\begin{tabular}{|c|l|c|c|}
\hline
Parameter & Description & Range & N\\
\hline
$\alpha$ & incident angle & $[\pi/12, \pi/3]$ & 30\\
$m_{G1}$ & frequency 1 in interface & $[10, 15]$ & 20\\
$m_{G2}$ & ratio of amplitudes of sinusoids & $[.5, 1]$ & 2\\
$m_{G3}$ & frequency 2 in interface & $[25, 30]$ & 20\\
\hline
 \end{tabular}
\end{center}

\end{table}

\section{Automatic classification}\label{sec:inverseproblem}

Define the forward prediction operator $F[{\bf m}]:X\rightarrow \RR^n$,
\begin{align}
F[{\bf m}]= \hat{Y}=(\hat{b}_{m^1}, \hdots, \hat{b}_{m^N}),\labeleq{prediction}
\end{align}
where the acoustic templates $b_{m^{i}}$ in the library $\mathcal{L}$ are given by \refeq{ym}. In material classification and geometric parameter estimation, the inverse problem is: From the backscattered signal  $Y\in \RR^{n}$ given by \refeq{backscattered} in the full seafloor domain, determine the seafloor parameters ${\bf m} = (m^{1}, \hdots, m^{N})$ with 
\begin{align}
F[{\bf m}]\simeq Y+\xi, \labeleq{invopt}
\end{align}
where $\xi$ is a random variable representing added noise.

Our optimization approach for $\refeq{invopt}$ involves first dividing the given data according to a partition of the full seafloor domain into subdomain of of width $\Delta s$. Moving from left to right across the full domain, each segment of the signal is matched with a signal in the library of acoustic templates. In determining a best fit, the algorithm relies on a data mistfit term involving the multilevel wavelet decomposition of the signals. 

\begin{figure}\caption{ The plot below shows the original signal $s$ and the coefficients in the multilevel wavelet decomposition using Haar wavelets. The original signal has length $2^{9}$, and the $l^{\text{th}}$ stage coefficients have length $2^{9-l}$, $1\leq l \leq 5$.}
\begin{center}

 \includegraphics[scale=.9]{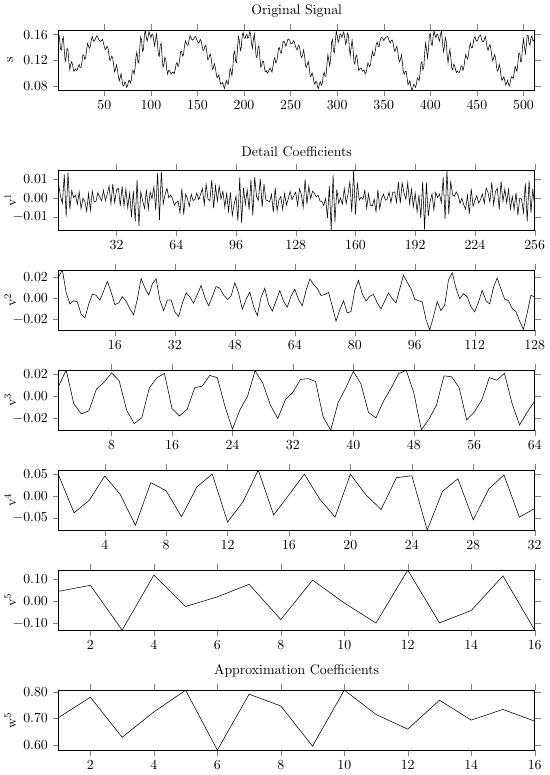}
\label{fig:wavelet-decomposition}
\end{center}
\end{figure}

Wavelet representations are used to extract information about signals according to the scale of variations within the signal \cite{Mallat1989}.  The multilevel discrete wavelet transform of a signal $s$ of length $2^{N_{1}}$ involves $N_{1}$ stages at most. In our experiments, we choose a maximum level $L_{max}<N_1$. In an application of the discrete wavelet transform of $s$, the signal is convolved with a lowpass filter and a highpass filter to produce vectors describing the seafloor at low frequencies and high frequencies respectively. These vectors are downsampled by a factor of 2, resulting in a vector of approximation coefficients $w^{1}$ and detail coefficients $v^{1}$, each of length $2^{N_{1}-1}$. The process is repeated for the vector of approximation coefficients, resulting in wavelet coefficients $w^{2}$ and $v^{2}$, each of size $2^{N_{1}-2}$. At the $l^{\text{th}}$ stage, the wavelet coefficients of $s$ are given by $w^{l}$ and $v^{l}$, each of size $2^{N_{1}-l}$. The final set of wavelet coefficients is a vector $w$ of length $2^{N_{1}}$ given by
\begin{align}
w=(w^{L_{max}}, v^{L_{max}}, v^{L_{max}-1}, \hdots, v^{1}).\labeleq{w}
\end{align}
The discrete wavelet transform preserves the $L_{2}-$ norm, that is, $\|s\|_{L_{2}}=\|w\|_{L_{2}}$.

Figure \ref{fig:wavelet-decomposition} contains plots of the original backscatter signal $s$, the approximation coefficient vector $w^{5}$, and individual plots of each vector $v^{i}$ of detail coefficients, $1\leq i \leq 5$. The $x$ axis corresponds to the coefficient index. 

\subsection{Algorithm}\label{sec:algorithm}
{The optimization process is a greedy table-lookup method, where regularization is imposed using constraints based on the parameters chosen for neighboring templates. Another option would be to use machine learning techniques such as neural networks in the identification \cite{Gorman1988}.} In each segment of the domain, the algorithm determines a match from the collection of templates in the library $\mathcal{L}$ for the truncated signal $(b_{i})_{j}=(Y)_{j}$ restricted to the $i^{th}$ segment of the measurement domain, $X\cap D^{i}$.
It is assumed that the observations $b_{i}$ and templates $\hat{b}_{m}\in \mathcal{L}$ are both vectors of length $2^{N_{1}}$ for an integer $N_{1}$. The wavelet coefficients of $b_{i}$ and $\hat{b}_{m}$ are given by $(w_i^l,v_i^l)$ and $(\hat{w}^l_{m},\hat{v}^l_{m})$, respectively.

The data misfit functions are defined for each segment of the domain and wavelet level,
\begin{align}
G^{l}_i(\hat{b}_{m}) = \|\hat{w}^{l}_{m}-w^{l}_{i}\|^{2}+\|\hat{v}^{l}_{m}-v^{l}_{i}\|^{2}.\labeleq{misfit}
\end{align}

Define $\mathcal{L}^{L_{max}+1}_{i}=\mathcal{L}$. The set $\mathcal{L}^{l}_{i}\subseteq{\mathcal{L}}^{l+1}_{i}$ is a collection of candidate templates $\hat{b}_{m}$ from the library that produce a sufficiently small value of the misfit function  given by \refeq{misfit}. For a given tolerance $\epsilon>0$, these sets are given by
\begin{align*}
\mathcal{L}^{l}_{i}=\{{\hat{b}_{m}}\in \mathcal{L}^{l+1}_i\mid G_i^l(\hat{b}_{m})<\epsilon2^{-l}\}.
\end{align*}

Since the maximum stage in the wavelet decomposition results in coefficient vectors of the smallest length, it is efficient to begin the matching process at the stage $l=L_{max}$. After the set $\mathcal{L}_{i}^{L_{max}}$ is determined, the decomposition level is decreased,  $l=L^{max}-1$,  and repeating the search process among the candidate templates $\hat{b}_{m}\in \mathcal{L}^{L_{max}}_{i}$, a refined set of candidate templates $\mathcal{L}^{L_{max}-1}_{i}\subseteq \mathcal{L}^{L_{max}}_{i}$ is produced.  This process is iterated for each further stage, $l=L_{max}-2, \hdots, l^{*}$, until either  $l^{*}=1$ or the pool of candidate templates $\mathcal{L}^{l^{*}-1}_i$ contains no templates.

If the set $\mathcal{L}^{l^{*}}_i$ contains multiple candidate templates, the choice for the seafloor parameter $m_{i}$ in the set $\mathcal{M}_{i}^{*}=\{m\in \mathcal{M}\mid \hat{b}_{m}\in \mathcal{L}^{l^{*}}_{i}\}$ is  made using information about spatially neighboring segments of the domain. There is no penalty term for the selection of ${m}^{1}$ corresponding to the leftmost segment of the domain. The penalty term  for $i>1$ is weighted by a fixed parameter $\delta>0$ and is given by
\begin{align*}
R_{i}(m', n')=\|m'_G-n'_{G}\|, 
\end{align*}
where $m'_{G}$ and $n'_{G}$ indicate the geometric parameters.

The algorithm for finding the best template match in each segment $i=1, \hdots N$, is summarized below.
\begin{enumerate}
\item Preliminary candidate search using multilevel wavelet decomposition: Set $l=L_{max}$.
\begin{enumerate}
\item  Candidate templates at level $l$ are given by  $\mathcal{L}^{l}_i$. 
\begin{enumerate}
\item If  $\mathcal{L}^{l}_i=\emptyset$, then $l^{*}=l+1$ and $\mathcal{L}_{i}^{{l^*}}$ is returned.
\item If $l=1$, then $l^{*}=1$ and $\mathcal{L}_{i}^{{l^*}}$ is returned. 
\end{enumerate}
\item  The wavelet level is decreased, $l=l-1$.
\end{enumerate}

\item The final parameter choice depends on the chosen parameters in the preceding segment of the domain,
\begin{align*}
{m}^{i}=\underset{m'\in\mathcal{M}^{*}_{i}}{\text{argmin}} \{G_{i}^{{l^*}}(\hat{b}_{m'})+ \delta R_{i}(m', m^{i-1})\}.
\end{align*}

\item The estimated seafloor parameters are ${\bf m}=(m^{1}, \hdots, m^{N})$
and the strength of the backscattered signal on the large domain is approximated by 
\begin{align*}
\hat{Y}=F[{\bf{m}}].
\end{align*}

\end{enumerate}

For a well-chosen tolerance level, the number of candidate templates will be significantly reduced in each iteration, resulting in the nested sets,
\begin{align*}
\mathcal{L}^{l^{*}}_i\subset \mathcal{L}^{l^{*}+1}_i\subset \hdots\subset \mathcal{L}^{L_{max}}_i \subset \mathcal{L}.
\end{align*}

The computational gain comes from this design; as the level of detail required to compare the signals is increased, the number of  possible candidate templates decreases.
The misfit function \refeq{misfit} is evaluated for all templates in the library only in the first iteration $l=L_{max}$, when the number of the wavelet coefficients in the vectors $\hat{w}^{l}_{m}$ and $\hat{v}^{l}_{m}$ is the smallest. As $l$ decreases, the number of wavelet coefficients in the vectors $\hat{w}^{l}_{m}$ and $\hat{v}^{l}_{m}$ increases, however the number of evaluations of the misfit \refeq{misfit} is lowered due to a reduced search space $\mathcal{L}^{l}_i\subset \mathcal{L}^{l+1}_i$.

\section{Numerical simulations}\label{sec:numerics}
The source frequency is 20 kHz and we set the incident angle of the incoming wave is $\alpha=\pi/6$.  The attenuation in the sediment layer for all of material types, listed in Table \ref{table:mat}, was set to $\gamma=$10 decibels per meter. From this, we determined that setting $h_{s}=2$ meters and $h_{w}=1$ meter resulting in a thin layer near the seafloor that sufficiently captured the high frequency oscillations and allowed for attenuation in the seafloor. The generic template domain size was set to $\Delta s=1$ meter to account for 20 wavelengths per template.

 Equations $\refeq{fwd1}-\refeq{fwd3}$ were discretized with linear finite elements using the COMSOL Multiphysics\textsuperscript{\textregistered} Acoustics Module \cite{Comsol}. Away from the interface, we used uniform meshes  with element size $\simeq \lambda/10$ in the water layer and $~\lambda/2$ in the sediment layer, where $\lambda$ is the minimum wavelength. We use a refined free triangulation for improved mesh quality  near the interface. The periodic boundary conditions on the left and right boundaries are modeled with Floquet periodic conditions in both the water and sediment domains. The computational domain is truncated at the top boundary using a $.1$ meter thick perfectly matched layer (PML). 

For the simulations of $\refeq{fwd1}-\refeq{fwd3}$ on the template domains, we found it numerically useful to enlarge the width of template domain to $4\Delta s$ meters, resulting in a complete mesh of 175481 elements and a linear system of equations with 356317 degrees of freedom. Then, the backscatter signals corresponding to the middle two segments, $\Delta s \leq x < 2\Delta s$ and $2\Delta s \leq x < 3\Delta s$ were each added as templates to the library corresponding to the same seafloor parameter. The enlarged domain accounts for numerical artifacts near boundary, and the overloading of templates adds redundancy to the library. MATLAB is used for the post-processing microlocal analysis step.

We tested the algorithm with backscatter signals generated from two seabeds of width of 20 meters, $D=[0 \text{m}, 20 \text{m}]\times [-2 \text{m}, 1 \text{m}]$. The varying material properties (sand, rock, and clay) and geometries in each model seabed were chosen to showcase the robustness of the proposed classification method:

\begin{itemize}
\item {\bf Model A.} This model has constant geometric parameters $m_{G1}\equiv 15 ,  m_{G2}\equiv 1$, $m_{G3}\equiv 26$, and material parameters
\begin{align}
m_{A}(x)=\begin{cases}\text{`sand'} & x\in[0, 2)\cup[8, 20] \\
\text{`rock'} & x\in[2, 5) \\
\text{`clay'} & x\in[5, 8)\end{cases}. \labeleq{modelA}
\end{align}

\item {\bf  Model B.} This model has smoothly varying geometric parameters $m_{G1}(x), m_{G2}(x)$, and $m_{G3}(x)$ that take values that do not exactly match the geometric parameters in the library. The material parameters are given by
\begin{align}
m_{A}(x)=\begin{cases}\text{`sand'} & x\in[0, 5) \\
\text{`rock'} & x\in[5, 13) \\
\text{`clay'} & x\in[13, 20)\end{cases}. \labeleq{modelB}
\end{align}

The geometric parameters are given by $m_{G1}(x)=14+ (.2 x-1)\chi_{[5,10)}(x)+ \chi_{[10,20)}(x)$, $
m_{G2}(x)=1+ (.5-.05 x)\chi_{[10,20)}(x)$, and 
$m_{G3}(x)=25+ x/20$.

\item \textbf{Buried objects.} Simulations of seafloors containing no buried objects as well as objects of varying width were simulated. The locations of the object domains added to the sediment domain in Model A and Model B are given by $D^{T}=[14-t_{w}, 14]\times [-h_s, -.02]$ and  $D^{T}=[4-t_{w}, 4]\times [-h_s, -.02]$, respectively. The object width $t_{w}$ takes values in the set $\{4, 2, 1, .5\}$, and the absence of the object implies $t_{w}=0$.

\end{itemize}
The synthetic backscattering is generated for both models from Helmholtz solutions $\refeq{fwd1}-\refeq{fwd3}$ on the full domain $D$. Then the response $d$ given by the formula \refeq{backscattered}, is polluted by 5\% Gaussian noise represented by the random variable $\xi$. The simulated backscatter data has the form
\begin{align}
b=d+\xi\labeleq{noisydata}
\end{align}

\noindent Then, the backscatter data was represented using a discrete wavelet representation with Haar wavelets and a maximum level $L_{max}=5$. The optimization procedure described in \S \ref{sec:inverseproblem} was implemented in MATLAB using the $L_{2}-$norm in with optimization parameters  $\delta=.02$ and $\epsilon=2^{-8}$.

\subsection{Lowering false alarm rate using an enriched library of templates}

We examined the effects of enriching the library using templates with multiple materials in terms of the rate of false detections of objects. In this experiment, we executed the classification algorithm on 100 different realizations of the data polluted with 5\% noise. To demonstrate the need for an enriched library, the algorithm was performed first on the library $\mathcal{L}_{0}$ containing templates with material parameters $m_{A}=$ `clay', `sand', `rock', and `metal'.  In Table \ref{table:falsealarm} the left column is a list of the transition types at locations $x=i\Delta s$ corresponding to the segments $D^{i}$ that were incorrectly identified by our method as `metal'. The second and third columns show the false alarm rate corresponding to Model A and Model B, respectively. The false alarm rate is defined as ratio of the number of times the segment  was labeled as a `metal' and the number of instances of that transition in the 100 trials. For example, there are instances of `sand-sand' in the material parameter \refeq{modelB} at the four locations $x=1, 2, 3, 4$. Out of 100 trials, there were three false detections at one of these locations, therefore the false alarm rate is $3/(100\times 4)$.

Then, the algorithm was performed on an enriched library $\mathcal{L}$ containing additional signals generated using synthetic seabeds containing a transition between two different materials, i.e., `sand-clay', `clay-sand', `clay-rock', `sand-metal'. The `sand' and `clay' templates generated from a seabed of length $2\Delta s$ containing the `sand-clay' and `clay-sand' transitions. Templates with material types `sand', `clay', and `metal' were added using seabeds containing the `clay-rock' and `sand-metal' type. Table \ref{table:falsealarm} indicates the main cause of false alarms is the presence of the `clay-rock' transition in the synthetic seabed. Segments of the domain that were a part of other transition types also contribute to the false alarm rate, however, the contribution is minimal in comparison. 

\begin{table}
\caption{Left: False alarm rate in the classification of seafloor parameters. The $'-'$ indicates that that transition was absent. First, the matching is performed using templates from a pure sediment library $\mathcal{L}_{0}$. Then, the matching is performed with an enriched library $\mathcal{L}$.  Right: The reconstruction algorithm is performed on 100 realizations of the signal \refeq{noisydata} corresponding to seabeds given in Model A and Model B. The object detection rate is the percentage of segments $D^{i}$ correctly labeled as `metal'}
\begin{center}
\begin{tabular}[t]{@{}*3l@{}}
\multicolumn{3}{c}{False Alarm Rate ($\mathcal{L}_{0}/\mathcal{L}$)} \\
 \toprule[0pt]
  \toprule[1pt]
\head{Transition type} &\head{Signal A} & \head{Signal B}  \\
 \cmidrule(r){1-1} \cmidrule(l){2-3}
sand-sand 		&0\%/0\% 	& 0.75\%/1.75\% 	\\
rock-rock  		& 0\%/0\% 		& 1.43\%/0.71\%	\\
sand-rock 		& 7\%	/0\% 	& -		\\
rock-clay 		& 10\%/0\% 	& -  		\\
clay-rock  		& -  		& 80\%/0\%        \\
\bottomrule[1pt]
\end{tabular}\hspace{10pt}
\begin{tabular}[t]{@{}*3l@{}}
\multicolumn{3}{c}{Detection Rate $(\mathcal{L}$)} \\
 \toprule[0pt]
  \toprule[1pt]
{\head{Object Width}} 
 &\head{Signal A} & \head{Signal B}  \\
 \cmidrule(r){1-1} \cmidrule(l){2-3}
4 meters & 99\% & 77.5\% \\
2 meters & 100\%  & 91\% \\
1 meter & 98\% & 73\% \\
0.5 meters & 1\% & 43\% \\
\bottomrule[1pt]
\end{tabular}

\end{center}
\label{table:falsealarm}
\end{table}%

\subsection{Inversion results}

Many well-studied approaches for determining the resolution needed to classify objects rely on ideas from information theory, for example, finding the resolution needed to to distinguish a sphere from a cube \cite{Kessel2002,  Myers2007,Pinto1997}. The proposed methods give a way to approach this problem using simulations of seabeds containing buried objects of varying size. Table \ref{table:falsealarm} shows the detection rate, defined as is the number of times a segment $D^{i}$ containing a submerged object was correctly identified divided by the number of trials.  

\begin{figure}

\caption{Shown in the plots below are synthetic backscatter data \refeq{noisydata} from seabeds given by Model A (left) and Model B (right) and submerged objects.  In red are the predictions  $F[\bf{m}]$ of the backscattering strength corresponding to the seafloor parameters ${\bf m}$ given by the reconstruction algorithm.} \label{fig:reconstruction}
\vspace{1em}
 \includegraphics[scale=1]{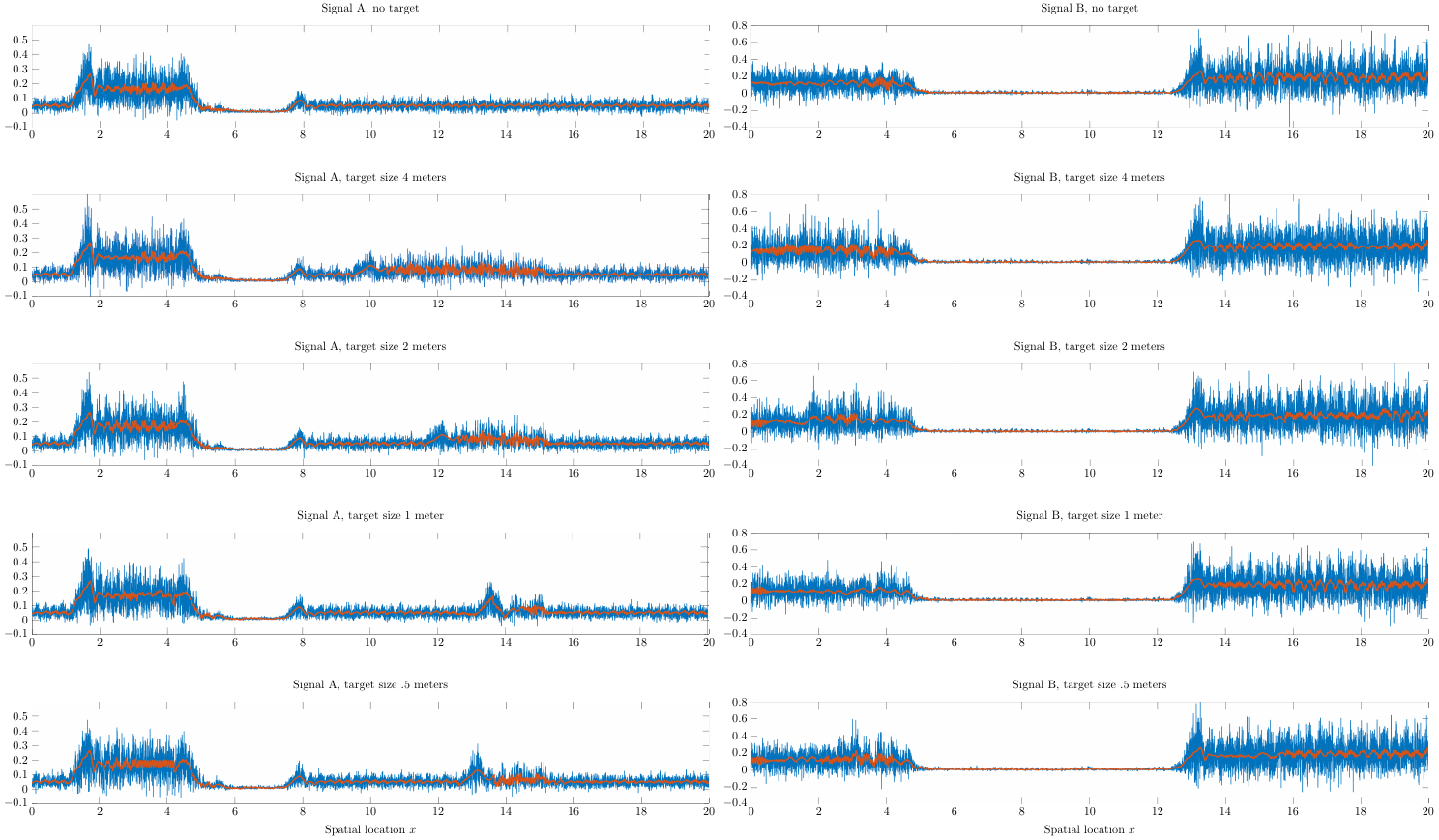}
  
\end{figure}
From these results, we conclude that the template domain size influences the recovery. For example, the smallest object with a width of .5 meters is rarely detected in Model A and not well-detected in Model B. In Model A, the objects can be detected by our algorithm to a resolution of $\Delta s=1$m with high accuracy. The higher detection rate for 2 meter objects over 4 meter objects in Model B might be explained by the location of the object near the boundary of the computational domain. In both cases, we believe that reducing the size of the template domain will increase the resolution of object detection. Also, we expect that enriching the library of templates to include more geometry parameters will improve the detection rate in Model B.  Figure \ref{fig:reconstruction} contains plots of the noisy signals \refeq{noisydata} and reconstructions using the proposed method. Though the geoacoustic properties of rock and metal are similar, the method correctly differentiated between the two.

%
%

\subsection{Reconstruction errors in seafloor geometry}
Though the primary aim of the proposed method is to properly identify objects based on material type, we also looked at how well the geometric parameters $m_{G}$ were recovered. In Model A, the geometric parameters are fixed and in Model B these parameters vary linearly across the domain. Figure \ref{fig:opterr} shows the distribution of relative errors in the reconstructed parameters $m_{G1}$, $m_{G2}$ and $m_{G3}$. The relative error formulas for the geometric parameters are
\begin{align*}
E_{i}=\frac{\|m_{Gi}-\hat{m}_{Gi}\|}{\sqrt{N}\|m_{Gi}\|},\qquad  i=1, 2, 3.
\end{align*}

In the library, there were two possibilities for the second geometric parameter, $m_{G2}=.5$ and $m_{G2}=1$. Since the seafloor model parameters \refeq{modelA} are constant, the error distribution of $E_{2}$ is split into correct estimates and incorrect estimates. The error distribution for Model B reflects the smooth changes in across the measurement domain. 

\begin{figure}\caption{Relative errors in the reconstruction of $m_{G1}$, $m_{G2}$, and $m_{G3}$.}

\begin{center}
 \includegraphics[scale=.8]{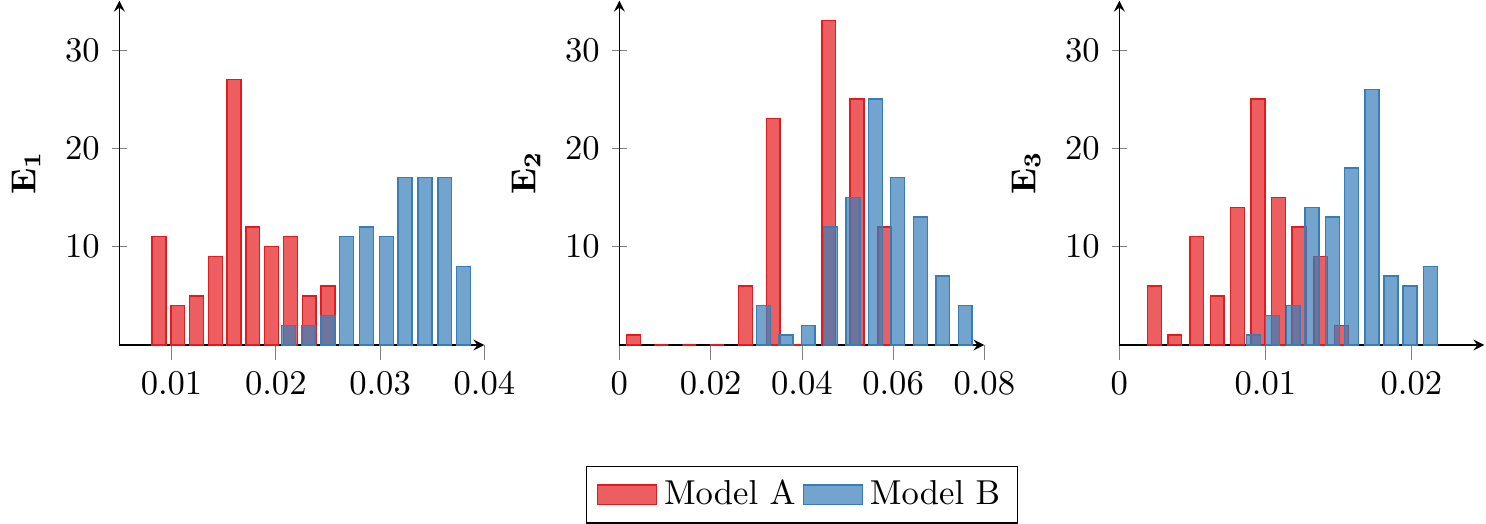}
\end{center}
\label{fig:opterr}

\end{figure}

\section{Conclusions and outlooks}\label{sec:conclusions}

The remote detection of {small features in the ocean} floor from measured acoustic backscatter remains a huge challenge in underwater acoustics. The proposed framework accomplishes this by matching data with a library of acoustic templates derived from multiscale simulations of acoustics. The efficiency is gained using a parametrization of local seafloor subdomains and using geometrical optics approximations away from the seafloor. The small size of the template domains is chosen to lower the computational cost of forward simulations, and a careful choice of parameterization results in fewer unknowns in the inverse problem. 

The templates are generated from forward simulations of Helmholtz equations corresponding to geometric and acoustic seafloor parameters corresponding to sand, clay, rock, and a submerged metal object. The range and resolution of the library was determined using realistic assumptions, requirements of the sonar instruments, and a parametric sensitivity study. In each patch of the seafloor, the matching algorithm uses the multilevel wavelet decomposition to determine good candidate templates. The process is enhanced using a smoothing term based on neighboring seafloor patches.

This work addresses problems that geometric optics cannot address at fine scales, for example, the detection of a submerged object. Here, only detailed simulations of Helmholtz equations can be used. Our results show that the material type of the sediment can be reliably determined by the proposed method. False alarms, or the incorrect identification of submerged metal objects, were lowered by using sediment templates that represented sediment layers consisting of a transition between two materials. We identified the clay-rock transition as the number one cause of false alarms. Templates representing other transitions were also used to enhance the classification. Numerical results indicate the accuracy of the reconstruction depends on the size of the template domains relative to the object size. 
%
\begin{figure}\caption{
{The large scale features in the synthetic image on the right comes from the material type and elevation (left, top layer) and the oscillations are due to the fine scale seafloor texture (left, bottom layer). The acoustic backscatter is simulated using an ensemble of solutions to $\refeq{fwd1}-\refeq{fwd3}$, with  $f(x,y, m_{G})=.01\sin(2\pi m_{G1}(x) x)+m_{G2}\sin(2\pi m_{G3}(x)x)$.} }\label{fig:2dresults}
\begin{center}
\resizebox {.8\columnwidth} {!} {
\includegraphics[height=5cm]{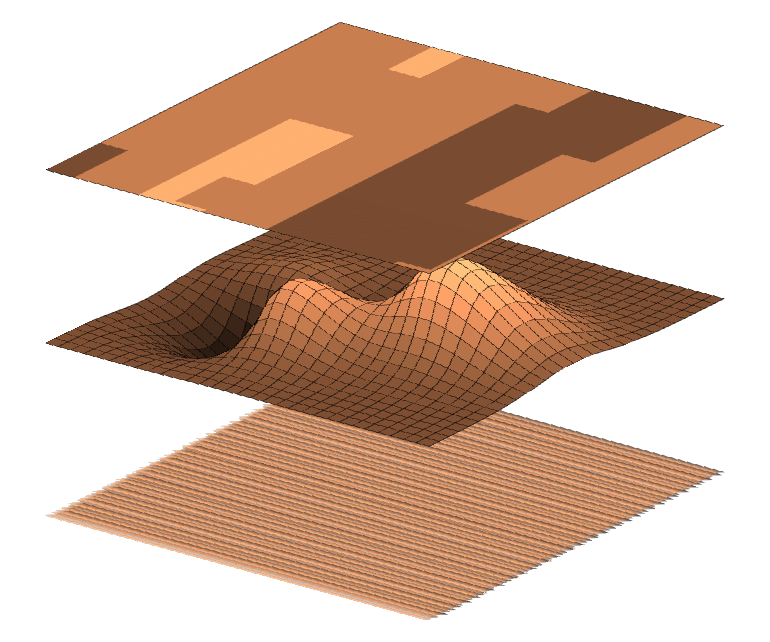}
\hspace{20pt}\includegraphics[width=5cm, height=5cm]{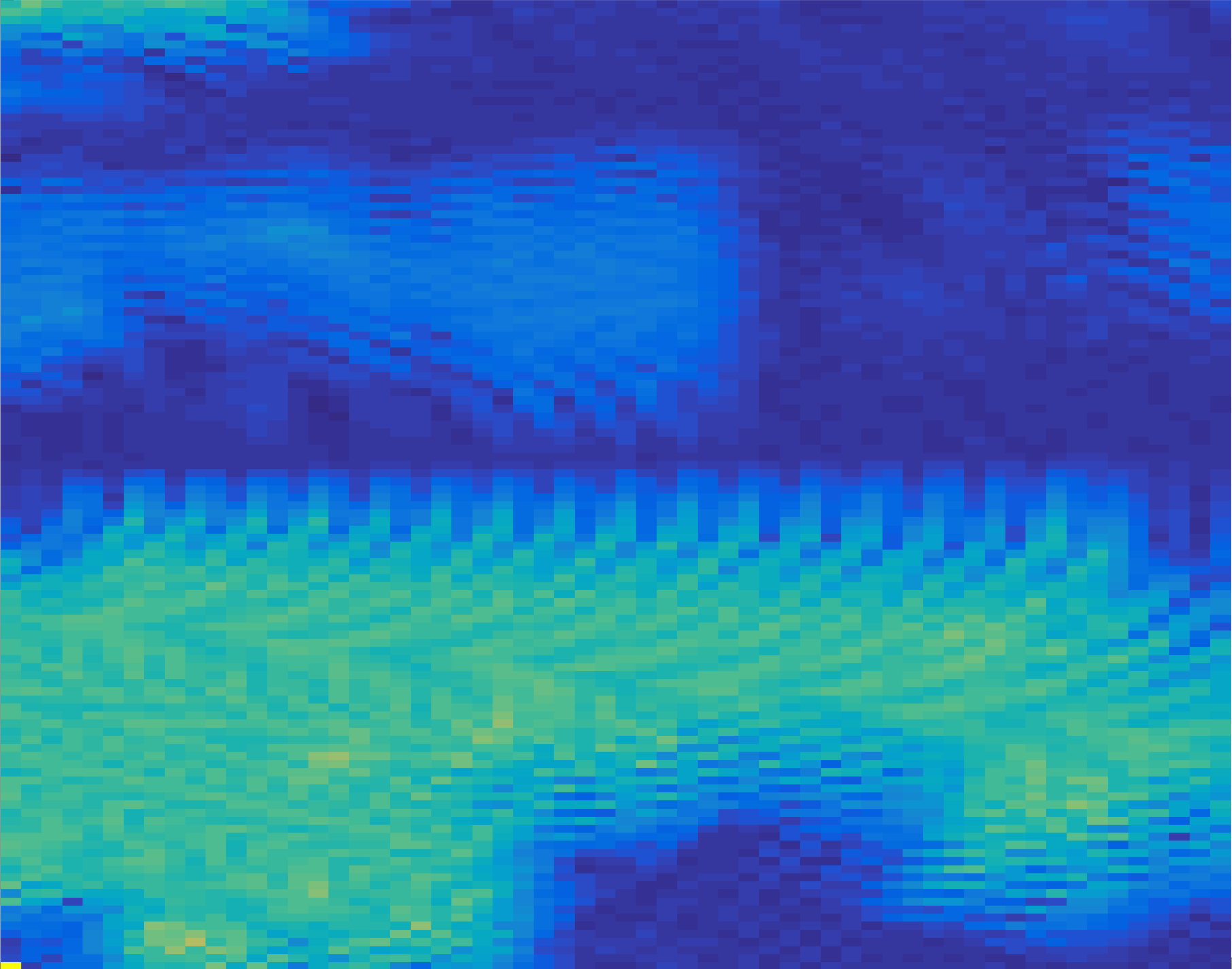}}
\end{center}
\end{figure}

The proposed framework can be extended using solutions Helmholtz equations in three dimensions. The forward model is based on sidescan sonars that illuminate a strip of the seafloor by transmitting an acoustic ping that acts as a spherical wave traveling toward the seafloor. The width of the strip depends on the incident angle $\alpha$, ranging from 15 degrees to 45 degrees, with smaller angles corresponding to higher frequencies.  In the model, the origin of the coordinate system is defined as be the projection of the source coordinate on the rough seafloor given by the function $z=f(x,y)$. Then, the synthetic backscatter datasets can be simulated using Helmholtz equations and numerical microlocal analysis in three dimensions, and an acoustic library can be formed using a parameterization of the sediment material type and the seafloor surface given by $z= f(x, y)$.  Simulations of these templates can be computed and stored efficiently, however, variable seabeds would introduce additional difficulties due to the complex scattering effects near transitions of sediment materials.

For demonstration, we approximate the process using simplified calculations, assuming that a reduced backscatter model is valid. Synthetic 2D images are generated using the existing techniques in the paper for determining the responses from regions of the sediment layer that are illuminated along the azimuth angles $\phi$, $0<\phi<\pi$ given by
\begin{align*}
z=f(r\cos(\phi), r\sin(\phi)), \qquad r_0\leq r \leq r_1,
\end{align*}
where the varying material type is now a function of both spatial coordinates $m_A=m_A(x, y)$.  Figure \ref{fig:2dresults} depicts a simulated 2D sonar image generated from a piecewise constant material parameter $m_{A}$.  The lower values (darker) represent fine sand or muddy sediment.

This study is a proof of concept of the inversion framework, and there are many possible extensions that are not in the scope of this paper but will be important to consider in the future. For more realistic problems, different techniques can be used. For example, the library of templates can be enriched using models that account for a wider range of sediments, objects with different materials and shapes, partially submerged objects, and templates adapted to different sonar system configurations. More parameters can be included in the inversion process, such as the attenuation, porosity, and properties of multiple horizontal sediment layers.

\section*{Acknowledgements}
We acknowledge the thoughtful suggestions and discussions with Gary Sammelmann, Frank Crosby, and Tory Cobb. The research was supported by the NSF grants DMS-1344199, DMS-1419027, and DMS-1620345,
and ONR Award N000141310408.

\section*{References}
\bibliography{multiscale-modeling-sonar}

\begin{thebibliography}{10}

\bibitem{Comsol}
COMSOL AB.
\newblock {COMSOL Multiphysics{\textregistered} v. 5.2. www.comsol.com.}

\bibitem{Baggeroer1988}
Arthur~B Baggeroer, W~A Kuperman, and Henrik Schmidt.
\newblock {Matched field processing: Source localization in correlated noise as
  an optimum parameter estimation problem}.
\newblock {\em Journal of Acoustic Society of America}, 83(2):571--587, 1988.

\bibitem{Bao1995}
Gang Bao.
\newblock Finite element approximation of time harmonic waves in periodic
  structures.
\newblock {\em {SIAM} Journal on Numerical Analysis}, 32(4):1155--1169, 1995.

\bibitem{Bao2005}
Gang Bao and Peijun Li.
\newblock Inverse medium scattering for the {Helmholtz} equation at fixed
  frequency.
\newblock {\em Inverse Problems}, 21:1621--1641, 2005.

\bibitem{Bao2013}
Gang Bao and Peijun Li.
\newblock Near-field imaging of infinite rough surfaces.
\newblock {\em SIAM Journal on Applied Mathematics}, 73(6):2162--2187, 2013.

\bibitem{Bell1997}
Judith~M. Bell.
\newblock {Application of optical ray tracing techniques to the simulation of
  sonar images}.
\newblock {\em Optical Engineering}, 36(6):1806--1813, 1997.

\bibitem{Benamou2004}
Jean~David Benamou, Francis Collino, and Olof Runborg.
\newblock {Numerical microlocal analysis of harmonic wavefields}.
\newblock {\em Journal of Computational Physics}, 199(2):717--741, 2004.

\bibitem{Reed1989}
T.~Beokett and Donald Hussong.
\newblock {Digital image processing techniques for enha , cement and
  classification of SeaMARC II side scan sonar imagery}.
\newblock {\em Journal of Geophysical Research}, 94(B6):7469--7490, 1989.

\bibitem{Blondel2009}
Philippe Blondel.
\newblock {\em {The Handbook of Sidescan Sonar}}.
\newblock Springer Berlin Heidelberg, Berlin, Heidelberg, 2010.

\bibitem{Bruno1991}
Oscar Bruno and Fernando Reitich.
\newblock {Solution of a boundary value problem for Helmholtz equation via
  variation of the boundary into the complex domain}.
\newblock 1991.

\bibitem{Burnett2015}
David~S Burnett.
\newblock {Computer simulation for predicting acoustic scattering from objects
  at the bottom of the ocean}.
\newblock {\em Acoustics Today}, 11(1):28--36, 2015.

\bibitem{Cervenka1993}
Pierre Cervenka and Christian {De Moustier}.
\newblock {Sidescan Sonar Image Processing Techniques}.
\newblock {\em IEEE Journal of Oceanic Engineering}, 18(2):108--122, apr 1993.

\bibitem{Chapman2001}
David M.~F. Chapman.
\newblock {What Are We Inverting For?}
\newblock In {\em Inverse Problems in Underwater Acoustics}, pages 1--14.
  Springer New York, New York, NY, 2001.

\bibitem{DelRioVera2009}
J~{Del Rio Vera}, E~Coiras, J~Groen, and B~Evans.
\newblock {Automatic Target Recognition in Synthetic Aperture Sonar Images
  Based on Geometrical Feature Extraction}.
\newblock {\em EURASIP Journal on Advances in Signal Processing}, 2009:1--10,
  2009.

\bibitem{Dobeck}
Gerald~J Dobeck.
\newblock {Algorithm Fusion for Automated Sea Mine Detection and
  Classification}.
\newblock In {\em MTS/IEEE Oceans 2001. An Ocean Odyssey. Conference
  Proceedings (IEEE Cat. No.01CH37295)}, pages 130--134, 2001.

\bibitem{Elston2004}
Gareth~R. Elston and Judith~M. Bell.
\newblock {Pseudospectral time-domain modeling of non-Rayleigh reverberation:
  Synthesis and statistical analysis of a sidescan sonar image of sand
  ripples}.
\newblock {\em IEEE Journal of Oceanic Engineering}, 29(2):317--329, apr 2004.

\bibitem{Fox1985}
Christopher~G. Fox and Dennis~E. Hayes.
\newblock {Quantitative methods for analyzing the roughness of the seafloor}.
\newblock {\em Reviews of Geophysics}, 23(1):1, 1985.

\bibitem{Frisk2009}
George~V. Frisk.
\newblock {Environmental influences on low‐frequency, shallow‐water
  acoustic propagation and inversion.}
\newblock {\em The Journal of the Acoustical Society of America},
  125(4):2590--2590, apr 2009.

\bibitem{Glassner1989}
Andrew Glassner, editor.
\newblock {\em {An Introduction to Ray Tracing}}, volume~8.
\newblock Morgan Kaufmann Publishers, San Francisco, 1989.

\bibitem{Gorman1988}
R.~Paul Gorman and Terrence~J. Sejnowski.
\newblock {Analysis of hidden units in a layered network trained to classify
  sonar targets}.
\newblock {\em Neural Networks}, 1(1):75--89, 1988.

\bibitem{Groen2010}
J.~Groen, E.~Coiras, J.~{Del Rio Vera}, and B.~Evans.
\newblock {Model-based sea mine classification with synthetic aperture sonar}.
\newblock {\em IET Radar, Sonar {\&} Navigation}, 4(1):62, 2010.

\bibitem{Hanrahan1993}
Pat Hanrahan and Wolfgang Krueger.
\newblock {Reflection from layered surfaces due to subsurface scattering}.
\newblock {\em Proceedings of the 20th annual conference on Computer graphics
  and interactive techniques}, pages 165--174, 1993.

\bibitem{Jackson2007}
Darrell~R. Jackson and Michael~D. Richardson.
\newblock {\em {High-Frequency Seafloor Acoustics}}.
\newblock Springer, New York, New York, USA, 2007.

\bibitem{Jackson1986}
Darrell~R. Jackson, D.P. Winebrenner, and A.~Ishimaru.
\newblock {Application of the composite roughness model to high-frequency
  bottom backscattering}.
\newblock {\em The Journal of the Acoustical Society of America}, 79(May):1410,
  1986.

\bibitem{Jensen2011}
Finn~B. Jensen, William~A. Kuperman, Michael~B. Porter, and Henrik Schmidt.
\newblock {Chapter 2 Wave Propagation Theory}.
\newblock In {\em Computational Ocean Acoustics}, pages 65--154. Springer New
  York, New York, NY, 2011.

\bibitem{Jensen2001}
Henrik~Wann Jensen, Stephen~R. Marschner, Marc Levoy, and Pat Hanrahan.
\newblock {A practical model for subsurface light transport}.
\newblock {\em Siggraph}, pages 511--518, 2001.

\bibitem{Kargl2014}
Steven~G. Kargl, Aubrey~L. Espana, Kevin~L. Williams, Jermaine~L. Kennedy, and
  Joseph~L. Lopes.
\newblock {Scattering From Objects at a Water???Sediment Interface: Experiment,
  High-Speed and High-Fidelity Models, and Physical Insight}.
\newblock {\em IEEE Journal of Oceanic Engineering}, 40(3):632--642, 2014.

\bibitem{Kargl2012}
Steven~G. Kargl, Kevin~L. Williams, and Aubrey~L. Espana.
\newblock {Fast model for target scattering in a homogeneous waveguide}.
\newblock {\em The Journal of the Acoustical Society of America}, 132(3):1909,
  2012.

\bibitem{Kessel2002}
Ronald~T. Kessel.
\newblock {Estimating the limitations that image resolution and contrast place
  on target recognition. Automatic Target Recognition XII}.
\newblock {\em Proceedings of SPIE}, 4726:316--327, jul 2002.

\bibitem{Landa2011}
Yanina Landa, Nicolay~M. Tanushev, and Richard Tsai.
\newblock {Discovery of point sources in the Helmholtz equation posed in
  unknown domains with obstacles}.
\newblock {\em Communications in Mathematical Sciences}, 9(3):903--928, 2011.

\bibitem{Lurton2010}
X.~Lurton.
\newblock {\em {An Introduction to Underwater Acoustics Principles and
  Applications}}.
\newblock Springer, 2010.

\bibitem{Mallat1989}
Stephane~G. Mallat.
\newblock {A Theory for Multiresolution Signal Decomposition: The Wavelet
  Representation}.
\newblock {\em IEEE Transactions on Pattern Analysis and Machine Intelligence},
  11(7):674--693, jul 1989.

\bibitem{Chapman1999}
Patrick~S. McCormick, Ji~Qiang, and Robert~D. Ryne.
\newblock {Visualizing underwater environments using multi-frequency sonar}.
\newblock {\em IEEE Computer Graphics and Applications}, 19(5):61--65, 1999.

\bibitem{Myers2007}
V.~Myers and M.~Pinto.
\newblock {Bounding the performance of sidescan sonar automatic target
  recognition algorithms using information theory}.
\newblock {\em IET Radar, Sonar {\&} Navigation}, 1(4):266--273, 2007.

\bibitem{Pailhas2010}
Yan Pailhas, Yvan Petillot, and Chris Capus.
\newblock {High-resolution sonars: What resolution do we need for target
  recognition?}, 2010.

\bibitem{Pinto1997}
M.~Pinto.
\newblock {Performance index for shadow classification in minehunting sonar -
  Google Scholar}.
\newblock In {\em Proceedings of the UDT Conference}, 1997.

\bibitem{Presig1994}
James~C. Preisig.
\newblock {A Minimax Approach to Adaptive Matched Field Processing in an
  Uncertain Propagation Environment}.
\newblock {\em IEEE Transactions on Signal Processing}, 42(6):1305--1316, 1994.

\bibitem{Reed2003}
Scott Reed, Yvan Petillot, and Judith Bell.
\newblock {An automatic approach to the detection and extraction of mine
  features in sidescan sonar}.
\newblock {\em IEEE Journal of Oceanic Engineering}, 28(1):90--105, jan 2003.

\bibitem{Rzhanov2012}
Yuri Rzhanov, Luciano Fonseca, and Larry Mayer.
\newblock {Construction of seafloor thematic maps from multibeam acoustic
  backscatter angular response data}.
\newblock {\em Computers and Geosciences}, 41:181--187, 2012.

\end{thebibliography}

\end{document}